\documentclass[12pt]{amsart}
\advance\hoffset by -0.5 truecm
\usepackage{amssymb}
\newtheorem{Theorem}{Theorem}[section]
\newtheorem{Lemma}[Theorem]{Lemma}

\newtheorem{Remark}[Theorem]{Remark}

\def\V{\mbox{Var}}

\def\Z{{\mathbb Z}}
\def\R\re
\def\V{\bf V}

\def \la{\lambda}

\def \re{{\mathbb R}}

\def \I{{\mathbb I}}

\def \0{\lambda_{0}}
\def \la{\lambda}
\def \ga{\gamma}

\def \G{{\mathbf G}_{M}}
\def\a{\omega_{1}}
\def\b{\omega_{2}}
\def\c{\omega_{3}}

\newcommand{\divv}{\overset v {\operatorname{div}}\,}
\newcommand{\divh}{\overset h {\operatorname{div}}\,}
\newcommand{\divm}{\overset m {\operatorname{div}}\,}
\newcommand{\de}[2]{\frac{\partial #1}{\partial #2}}
\newcommand{\del}[2]{\frac{\delta #1}{\delta #2}}
\newcommand{\ov}{\overset v \omega}
\newcommand{\oh}{\overset h \omega}

\begin{document}
\title[]{Rigidity properties of Anosov optical hypersurfaces}

\author[N.S. Dairbekov]{Nurlan S. Dairbekov}
\address{Kazakh British Technical University,
Tole bi 59, 050000 Almaty, Kazakhstan }
\email{Nurlan.Dairbekov@gmail.com}

\author[G.P. Paternain]{Gabriel P. Paternain}
 \address{ Department of Pure Mathematics and Mathematical Statistics,
University of Cambridge,
Cambridge CB3 0WB, England}
 \email {g.p.paternain@dpmms.cam.ac.uk}




\begin{abstract} We consider an optical hypersurface $\Sigma$ in
the cotangent bundle $\tau:T^*M\to M$ of a closed manifold $M$
endowed with a twisted symplectic structure. We show that if the
characteristic foliation of $\Sigma$ is Anosov, then a smooth
1-form $\theta$ on $M$ is exact if and only $\tau^*\theta$ has
zero integral over every closed characteristic of $\Sigma$. This
result is derived from a related theorem about magnetic flows
which generalizes our work in \cite{DP}. Other rigidity issues are
also discussed.

\end{abstract}

\maketitle

\section{Introduction}

Let $M$ be a closed connected $n$-manifold and let $\tau:T^*M\to M$ be its cotangent bundle.
Given an arbitrary smooth closed 2-form $\Omega$ on $M$, we consider $T^*M$ endowed with the
{\it twisted symplectic structure}
\[\omega:=-d\lambda+\tau^*\Omega,\]
where $\la$ is the Liouville 1-form.\footnote{Hence, we use the
convention that the Hamiltonian vector field $X_{H}$ of a
Hamiltonian $H$ is determined by $i_{X_{H}}\omega=dH$.}

A smooth, closed, connected, fiberwise strictly convex hypersurface $\Sigma \subset T^*X$ is
called \emph{optical}.\footnote{For the origins of the term {\it optical} see \cite[Section 9]{Arn}.}
Fiberwise strict convexity means that $\Sigma$ intersects
each fiber $T_x^*X$ along a hypersurface whose second fundamental form is positive
definite. Denote by $\sigma$ the characteristic foliation of $\Sigma$, i.e., the
1--dimensional foliation tangent to the kernel of $\omega|_{T\Sigma}$. Note that
$\sigma$ is orientable.

We shall say that an optical hypersurface $\Sigma\subset T^*M$ is {\it Anosov} (or hyperbolic)
if the characteristic foliation admits a (non-vanishing) tangent vector field whose flow
is Anosov. Since the flows of two such vector fields are reparametrizations of one another, the property
of being Anosov is independent of the chosen vector field (cf. \cite{AS}) and is a property of $\Sigma$.

In the present paper we shall study various rigidity properties of Anosov optical hypersurfaces
on cotangent bundles equipped with twisted symplectic structures. These properties
are motivated by recent results that we obtained for two dimensional magnetic flows \cite{DP}.

Here is one of our main results:

\medskip

\noindent {\bf Theorem A.} {\it Let $\Sigma\subset T^*M$ be an Anosov optical hypersurface, where
$T^*M$ is endowed with a twisted symplectic structure $-d\la+\tau^*\Omega$.
Let $\theta$ be a smooth 1-form on $M$. Then $\theta$ is exact if and only if
\[\int_{\Gamma}\tau^*\theta=0\]
for every closed characteristic $\Gamma$ of $\sigma$.
}

\medskip

If the closed 2-form $\Omega$ determines an integral class, we can
introduce the notion of action spectrum as follows. Suppose
$[\Omega]\in H^2(M,\Z)$. Then there exists a principal circle
bundle $\Pi:P\to M$ with Euler class $[\Omega]$. The bundle admits
a connection 1-form $\alpha$ such that $d\alpha=-2\pi
\,\Pi^*\Omega$. Let $\log {\mbox{\rm hol}}_{\alpha}:Z_{1}(M)\to
\re/\Z$, be the logarithm of the holonomy of the connection
$\alpha$. Here, $Z_1(M)$ is the space of 1-cycles and for every
2-chain $f:\Sigma\to M$ we have
\[\log {\mbox{\rm hol}}_{\alpha}(\partial\Sigma)=-\,\int_{\Sigma}f^*\Omega\;\;{\mbox{\rm mod}}\,1.\]

We define the {\it action} of an oriented closed characteristic $\Gamma$ as:
\[A(\Gamma):=\int_{\Gamma}\la+\,\log {\mbox{\rm hol}}_{\alpha}(\tau(\Gamma))\;{\mbox{\rm mod}}\,1.\]
We call the set ${\mathcal S}\subset \re/\Z$ of values $A(\Gamma)$ as $\Gamma$
ranges over all (oriented) closed characteristics, the {\it action spectrum} of $\Sigma$.

If $\Omega$ does not determine an integral class, but there exists $c\neq 0$ such that
$[c\Omega]\in H^{2}(M,\Z)$ we can still define the action spectrum by considering $R_{c}(\Sigma)$
and $-d\la+c\tau^*\Omega$, where $R_{c}(x,p):=(x,cp)$. The characteristic foliations of
$(\Sigma, -d\la+\tau^*\Omega)$ and $(R_{c}(\Sigma), -d\la+c\tau^*\Omega)$ are conjugate
by $R_{c}$.

Suppose now that we vary
the connection 1-form $\alpha$. Let $\alpha_r$ be a smooth 1-parameter
family of connections for $r\in (-\varepsilon, \varepsilon)$ with
$\alpha_0=\alpha$
Then we can write $\alpha_{r}-\alpha=\Pi^*\beta_r$, where $\beta_r$
are smooth 1-forms on $M$.
The connection $\alpha_r$ has curvature form $-2\pi \,\Omega+d\beta_r$.
If we let $\Omega_{r}=\Omega-\frac{1}{2\pi }d\beta_r$
we get a characteristic foliation $\sigma^r$ and an action spectrum
${\mathcal S}_{r}$. If the characteristic foliation $\sigma$ is Anosov, then for
$\varepsilon$ small enough $\sigma^r$ is Anosov for all
$r\in (-\varepsilon, \varepsilon)$.

\medskip

\noindent {\bf Corollary 1.} {\it Let $M$ be a closed connected manifold
and let $\Sigma\subset T^*M$ be an optical hypersurface.
Let $\Omega$ be a closed integral $2$-form and suppose that
$(\Sigma, -d\la+\tau^*\Omega)$ is Anosov.
If ${\mathcal S}_{r}={\mathcal S}$ for all
$r$ sufficiently small, then the deformation is trivial, that is,
$\alpha_r=\alpha+\Pi^{*}dF_r$ and $\Omega_r=\Omega$, where $F_{r}$ are smooth
functions on $M$.}

\medskip

The proof of Corollary 1 is very similar to that of Theorem C in \cite{DP} and hence we omit it.

Theorem A will be a consequence of the following result. Let $M$ be a closed connected manifold
endowed with a Finsler metric $F$. The Legendre transform $\ell_{F}:TM\setminus\{0\}\to T^*M\setminus\{0\}$
associated with the Lagrangian $\frac{1}{2}F^2$ is a diffeomorphism and
$\omega_{0}:=\ell_{F}^{*}(-d\la)$ defines a symplectic form on $TM\setminus\{0\}$. Now let $\Omega$
be a smooth closed 2-form on $M$ and $\pi:TM\to M$ the canonical projection.
The {\it magnetic flow} of the pair $(F,\Omega)$ is the Hamiltonian flow $\phi$ of
$\frac{1}{2}F^2$ with respect to the symplectic form $\omega_{0}+\pi^*\Omega$.
We shall consider $\phi$ restricted to the unit sphere bundle $SM:=F^{-1}(1)$.
A curve $\ga:\re\to M$ given by $\ga(t)=\pi(\phi_{t}(x,v))$
will be called a {\it magnetic geodesic}.

\medskip

\noindent {\bf Theorem B.} {\it Let $(M,F)$ be a closed connected Finsler manifold and $\Omega$
an arbitrary smooth closed $2$-form. Suppose the magnetic flow $\phi$ of the pair
$(F,\Omega)$ is Anosov and let $\G$ be the vector field generating $\phi$.

If $h:M\to\mathbb R$ is any smooth function and $\theta$ is any smooth $1$-form on $M$ such that
there is a smooth function $u:SM\to \mathbb R$ for which
$h(x)+\theta_x(v)=\G(u)$, then
$h$ is identically zero and $\theta$ is exact.

}

\medskip

Note that by the smooth Liv\v sic theorem \cite{LMM} saying that
$h(x)+\theta_x(v)=\G(u)$ is equivalent to saying that
$h(x)+\theta_x(v)$ has zero integral over every closed magnetic geodesic.

Various versions of Theorem B were previously known:

\begin{enumerate}
\item V. Guillemin and D. Kazhdan in \cite{GK1} proved Theorem B for $M$ a surface, $\Omega=0$
and $F$ a negatively curved Riemannian metric. In \cite{GK2} they extended this to higher dimensional manifolds
under a pointwise curvature pinching assumption and Min-Oo \cite{min} proved it when the curvature operator is
negative definite. All these results were based on Fourier analysis.
\item A major breakthrough was obtained by C. Croke and V. Sharafutdinov \cite{CS} in which results like
Theorem B were proved just assuming negative sectional curvature and in any dimension. The novel ingredient
here was the {\it Pestov identity}.
\item In \cite{DS}, Dairbekov and Sharafutdinov proved Theorem B, just assuming that the geodesic flow
of the Riemannian metric is Anosov.
\item In \cite{DP}, the authors proved Theorem B when $M$ is a surface and $F$ is a Riemannian metric, but
$\Omega$ is arbitrary.
\end{enumerate}

We now describe some applications of these results.

\subsection{Infinitesimal spectral rigidity}

Corollary 1 and the results of V. Guillemin and A. Uribe in \cite{GU} give a version of infinitesimal
spectral rigidity for magnetic flows. This version was obtained in \cite{DP}
for the case of surfaces. Suppose $\Omega$ is a closed integral 2-form and $g$
a Riemannian metric.
For every positive integer $m$, let $L_m$ be the Hermitian line bundle with connection
over $M$ associated with $\Pi$ via the character $e^{i\theta}\mapsto e^{i m\theta}$ of $S^1$.
The metric on $M$, together with the connection on $L_m$ determine a Bochner-Laplace operator
acting on sections of $L_m$. For each $m$, let $\{\nu_{m,j}:\,j=1,2,\dots\}$ be the spectrum
of this operator. If we now vary the connection 1-form $\alpha$ as above we obtain
eigenvalues $\nu_{m,j}^{r}$.

\medskip

\noindent {\bf Corollary 2.} {\it Let $M$ be a closed connected manifold
endowed with a Riemannian metric $g$ and let $\Omega$ be an integral $2$-form.
Suppose the magnetic flow of the pair $(g,\Omega)$ is Anosov.
If $\nu_{m,j}^r$ is independent of $r$ for all $m$ and $j$ $($i.e. the deformation
is isospectral$)$, then the deformation is trivial, that is,
$\alpha_r=\alpha+\Pi^{*}dF_r$ and $\Omega_r=\Omega$, where $F_{r}$ are smooth
functions on $M$.}

\medskip

Indeed, let us consider the periodic distribution
\[\Upsilon(s)=\sum_{m,j}\varphi\left(\sqrt{\nu_{m,j}+m^2}-m\sqrt{2}\right)e^{ims}\]
where $\varphi$ is a Schwartz function on the real line.
Theorem 6.9 in \cite{GU} asserts that the singularities of $\Upsilon$ are included in the
set of all $s\in\re$ for which $s/2\pi\,\mbox{\rm mod}\,1\in {\mathcal S}$. Moreover, each point
of the action spectrum arises as a singularity of $\Upsilon$ for an appropriate choice of $\varphi$.
Corollary 2 is now an immediate consequence of Corollary 1.

There is an equivalent way of formulating Corollary 2 in purely
Riemannian terms using the {\it Kaluza-Klein metric}. Consider on
$P$ the metric $g_{KK}$ defined uniquely by the following
conditions: the restriction of $d\Pi$ to the horizontal subspace
of the connection $\alpha$ is an isometry, vertical and horizontal
subspaces are orthogonal and the vector field
$\partial/\partial\theta$ tangent to the fibres has norm one. If
we vary the connection $\alpha$ as above we obtain a 1-parameter
family of Kaluza-Klein metrics $g_{KK}^r$, $r\in
(-\varepsilon,\varepsilon)$. Consider the usual Laplacian
$\Delta_{KK}^r$ of these metrics. Corollary 2 could be rephrased
by saying that if the spectrum of $\Delta_{KK}^r$ remains
unchanged, then the deformation is trivial. In fact, the
eigenvalues $\lambda_{m,j}$ of $\Delta_{KK}$ restricted to the
$(-m)$-eigenspace of $-i\partial/\partial\theta$ are related to
$\nu_{m,j}$ by $\lambda_{m,j}=\nu_{m,j}+m^2$, cf. \cite[Section
6]{GU}.

\bigskip

\subsection{Regularity of the Anosov splitting} Theorem B can be used for the study
of the regularity of the Anosov splitting of magnetic flows. In
fact, in dimension two this problem is completely solved in the
Riemannian setting in \cite{DP} and is one of the main motivations
of this paper. Here we show:

\medskip

\noindent {\bf Theorem C.} {\it Let $M$ be a closed connected
manifold endowed with a Finsler metric $F$ and let $\Omega$ be an exact $2$-form.
Suppose that the magnetic flow $\phi$ of the pair $(F,\Omega)$ is Anosov.
If the Anosov splitting of $\phi$ is of class $C^1$, then $\Omega$ must vanish, i.e.,
the magnetic flow is a Finsler geodesic flow.
}

\medskip

Theorem C was proved in \cite{Pa}, when $F$ is a Riemannian metric, using Aubry-Mather theory.
The proof in \cite{Pa} {\it cannot} be extended to include arbitrary (non-reversible) Finsler
metrics, since it uses the invariance of the Riemannian metric under the flip $(x,v)\mapsto (x,-v)$.

\subsection{Sketch of the proof of Theorem B}

Perhaps the most important element in the proof is the Pestov identity in our setting.
This comes in two flavours. We first obtain a scalar identity (cf. Lemma \ref{pes-2} in dimension
two and Lemma \ref{pes-n} in arbitrary dimension). When this identity is manipulated and integrated with
respect to the Liouville measure $\mu$ of $SM$ it gives rise to our key integral identity:

\begin{multline}\label{pestov-integral-final}
\int_{SM}\big\{|\mathbf X (\nabla^{\cdot} u)|^2
-\langle \mathbf R_y(\nabla^{\cdot} u),\nabla^{\cdot} u\rangle
-L(Y(y),\nabla^{\cdot} u,\nabla^{\cdot}u)
-\langle\nabla^{\cdot } (\mathbf X u),Y(\nabla^{\cdot} u)\rangle
\\
-2\langle Y(y),\nabla^{\cdot}u\rangle^2
+\langle \nabla^{:} u,Y(\nabla^{\cdot}u)\rangle
+\langle\nabla_{|(\nabla^{\cdot}u)}Y(y),\nabla^{\cdot}u\rangle\big\}\,d\mu
\\
=\int_{SM}\big\{|\nabla^{\cdot} (\mathbf X  u)|^2
-n(\mathbf X u)^2\big\}\,d\mu.
\end{multline}
Of course, this formula needs explaining and we shall fully do so in Sections 3 and 4, but for the purpose of this
sketch it suffices to note the following points:
\begin{enumerate}
\item $u$ is a smooth positively homogeneous function of degree
zero on $TM\setminus\{0\}$; \item $\mathbf X$ is a suitable vector
field on $TM\setminus\{0\}$ whose restriction to $SM$ coincides
with $\G$; \item the various derivatives that appear in the
formula are all obtained using the Chern connection of the Finsler
metric and are explained in detail in Section 4; \item
$\nabla^{\cdot}u$ vanishes if and only if $u$ is the pull back of
a function on $M$; \item inner products and norms are all taken
with respect to the fundamental tensor in Finsler geometry:
$$
g_{ij}(x,y)=\frac12[F^2]_{y^iy^j}(x,y);
$$
\item ${\bf R}$ and $L$ are respectively the Riemann curvature operator and the
Landsberg tensor from Finsler geometry; $Y$
is the Lorentz force associated with the magnetic field;
\item $n$ is the dimension of $M$.
\end{enumerate}

We may regard the identity as a kind of ``dynamical Weitzenb\"ock
formula". Suppose now that $\G(u)=h\circ\pi+\theta$ and extend $u$
to a positively homogeneous function of degree zero on
$TM\setminus\{0\}$ (still denoted by $u$).  Then $\mathbf X(u)=F
h\circ\pi+\theta$ and it is not hard to see (cf. Lemma
\ref{int-nabla}) that the right hand side of
(\ref{pestov-integral-final}) is non-positive and thus

\begin{multline*}
\int_{SM}\big\{|\mathbf X (\nabla^{\cdot} u)|^2
-\langle \mathbf R_y(\nabla^{\cdot} u),\nabla^{\cdot} u\rangle
-L(Y(y),\nabla^{\cdot} u,\nabla^{\cdot}u)
-\langle\nabla^{\cdot } (\mathbf X u),Y(\nabla^{\cdot} u)\rangle
\\
-2\langle Y(y),\nabla^{\cdot}u\rangle^2
+\langle \nabla^{:} u,Y(\nabla^{\cdot}u)\rangle
+\langle\nabla_{|(\nabla^{\cdot}u)}Y(y),\nabla^{\cdot}u\rangle\big\}\,d\mu
\leq 0.
\end{multline*}

It is at this point that we need a new ingredient. We will note
that the left hand side of the last inequality is closely related
to an analogue of the classical index form in Riemannian geometry.
Bilinear forms of this type already appeared in \cite{PP1} and
were very useful for the study of derivatives of topological
entropy. This time the form that we need is a sharper version of
the one that appears in \cite{PP1}. The key point is that the
Anosov property, via the abscence of conjugate points established
in \cite{PP,P3} (see \cite{CGIP} for a proof using the asymptotic
Maslov index), will imply that when we integrate the expression
inside the brackets in the last inequality along every closed
magnetic geodesic the outcome should be {\it non-negative} and
{\it zero} if and only if $\nabla^{\cdot} u$ vanishes along every
closed magnetic geodesic. When we combine this fact with the
recent {\it non-negative Liv\v sic theorem} \cite{LT,PS} we deduce
that $\nabla^{\cdot} u$ must vanish over every closed magnetic
geodesic and thus it must be identically zero on
$TM\setminus\{0\}$. This means that $u=f\circ\pi$ where $f$ is a
smooth function on $M$. But in this case, since
$d\pi_{(x,v)}(\G)=v$ we have $\G(u)=df_{x}(v)$ and Theorem B
follows.

A considerable part of the paper will be devoted to the proof of
the integral formula (\ref{pestov-integral-final}). This necessitates the language and formalism
of Finsler geometry which makes the derivation of the formula a bit cumbersome.
To help the reader, we have included a brief section in which we prove the integral
formula in dimension two. This easier case still shows some of the main features
and it can be read independently of the other sections.


\section{Theorem B implies Theorem A} Let us explain why Theorem B implies Theorem A.

Suppose $\Sigma\subset T^*M$ is an optical hypersurface which encloses an open region $U$ in
$T^*M$. Let $\Sigma_{x}:=\Sigma\cap T^{*}_{x}M$ which
is a strictly convex hypersurface in the vector space $T^{*}_{x}M$ which encloses
$U_{x}:=U\cap T_{x}^*M$. Consider an auxiliar smooth Riemannian
metric $g$ on $\tau:T^*M\to M$, that is, for each $x\in M$, $g_{x}$ is an inner product in $T_{x}^{*}M$.
For each $x\in M$, the inner product $g_x$ gives rise to a volume form $\varpi_{x}$ in $T_{x}^*M$.
Consider the barycenter of $U_{x}$, i.e.,
\[\beta_{x}:=\frac{\int_{U_{x}}p\,\varpi_{x}}{\int_{U_{x}}\varpi_{x}}.\]
The map $x\mapsto\beta_{x}$ can be seen as a smooth 1-form and by strict convexity $\beta_{x}\in U_{x}$
for all $x\in M$.

Consider the map $B:T^*M\to T^*M$ given by $B(x,p)=(x,p-\beta_{x})$. It is easy to check that
$B^*(\la)=\la-\tau^*\beta$ and that $B^*(\tau^*\Omega)=\tau^*\Omega$.
Hence if we let $\widetilde{\Omega}:=\Omega+d\beta$, $B$ is
a symplectomorphism between $(T^*M,-d\la+\tau^*\Omega)$ and $(T^*M,-d\la+\tau^*\widetilde{\Omega})$.
Now set $\widetilde{\Sigma}:=B(\Sigma)$ and observe that $\widetilde{\Sigma}$ is optical and contains the zero
section of $T^*M$. Also note that
\[\int_{\Gamma}\tau^*\theta=0\]
for all $\Gamma$ of $\sigma$ if and only if
\[\int_{\widetilde{\Gamma}}\tau^*\theta=0\]
for all $\widetilde{\Gamma}$ of $\widetilde{\sigma}$.
Thus, without loss of generality, we may assume that $\Sigma$ contains the zero section of
$T^*M$. But in that case we can define a Finsler metric $F$ on $M$
using homogeneity and declaring that $\Sigma$ corresponds to the unit cosphere bundle of $F$.
The hypothesis in Theorem A tells us that
\[\int_{\gamma}\theta=0\]
for every closed magnetic geodesic $\ga$ of $(F,\Omega)$. The smooth Liv\v sic theorem \cite{LMM} and Theorem B
imply that $\theta$ must be exact.

\section{Proof of Theorem B for surfaces}

\subsection{Canonical coframing}

Let $M$ be a closed oriented connected surface.
A smooth {\it Finsler structure} on $M$ is
 a smooth hypersurface $SM\subset TM$
for which the canonical projection $\pi:SM\to M$ is a surjective
submersion having the property that for each $x\in M$, the $\pi$-fibre
 $\pi^{-1}(x)=SM\cap T_{x}M$ is a smooth, closed, strictly
 convex curve enclosing
the origin $0_x\in T_{x}M$.

Given such a structure it is possible to define a canonical coframing
$(\a,\b,\c)$ on $SM$ that satisfies the following structural
 equations (see \cite[Chapter 4]{BCZ}):
\begin{align}
d\a&=-\b\wedge\c,\label{eq1}\\
d\b&=-\c\wedge(\a-I\b),\label{eq2}\\
d\c&=-(K\a-J\c)\wedge\b.\label{eq3}
\end{align}
where $I$, $K$ and $J$ are smooth functions on $SM$.
The function $I$ is called the {\it main scalar}
of the structure. When the Finsler structure is Riemannain,
$K$ is the Gaussian curvature.

The form $\a$ is the canonical contact form of $SM$ whose Reeb vector
field is the geodesic vector field $X$.
The volume form $\a\wedge d\a$ gives rise to the Liouville measure
$d\mu$ of $SM$.

Consider the vector fields $(X,H,V)$ dual to
$(\a,\b,\c)$. As a consequence of (\ref{eq1}--\ref{eq3}) they satisfy the commutation relations
\begin{equation}\label{comm}
[V,X]=H,\quad [H,V]=X+IH+JV,\quad [X,H]=KV.
\end{equation}

Below we will use the following general fact. Let $N$ be a closed
oriented manifold and $\Theta$ a volume form. Let $X$ be a vector
field on $N$ and $f:N\to\re$ a smooth function. Then
\begin{equation}
\int_{N}X(f)\,\Theta=-\int_{N}f\,L_{X}\Theta,
\label{bas}
\end{equation}
where $L_X\Theta$ is the Lie derivative of $\Theta$ along $X$.

Now let $\Theta:=\a\wedge\b\wedge\c$. Using the commutation relations we obtain:

\begin{align}
L_{X}\Theta &=0;\label{lie1}\\
L_{H}\Theta &=-J\,\Theta;\label{lie2}\\
L_{V}\Theta &=I\,\Theta.\label{lie3}
\end{align}

\subsection{Identities} Let $\Omega$ be a 2-form on $M$. An important observation is this:
$\pi^*\Omega=\la\,\omega_1\wedge\omega_2$, where $\la:SM\to\re$ is a function such that
\[V(\la)=-\la\,I.\]
This relation is obtained using the structure equations in
$d(\la\,\omega_1\wedge\omega_2)=0$. The magnetic vector field is
$$\G=X+\la V.$$
The brackets are now:
\begin{equation}
[V,\G]=H-\la IV,\quad [H,V]=\G+IH+(J-\la)V,\quad [\G,H]={\mathbb K}V-\la \G-\la IH,
\end{equation}
where ${\mathbb K}:=K-H(\la)+\la^{2}-\la J$.

Using these brackets we obtain as in \cite[Lemma 3.1]{DP}:

\begin{Lemma}[The Pestov identity]\label{pes-2}
For every smooth function $u:SM\to\mathbb R$ we have
\begin{align*}
2Hu\cdot V\G u
&=(\G u)^2+(Hu)^2-{\mathbb K}(Vu)^2
+\G(Hu\cdot Vu)\\
&-H(\G u\cdot Vu)
+V(\G u\cdot Hu)+\G u\cdot(IHu+JVu).
\end{align*}
\end{Lemma}

We omit the proof which is (once you know the formula!) a straightforward verification
using the bracket relations.

Integrating Pestov's identity over $SM$ against the Liouville
measure $d\mu$ and using (\ref{bas}) and (\ref{lie1}--\ref{lie3})
we obtain:
\begin{equation}\label{id1}
2\int_{SM} Hu\cdot V\G u\,d\mu
=\int_{SM}(\G u)^2\,d\mu+\int_{SM}(Hu)^2\,d\mu
-\int_{SM}{\mathbb K}(Vu)^2\,d\mu.
\end{equation}

By the commutation relations, we have
$$
\G Vu=V\G u-Hu+\la IVu.
$$
Therefore,
\[(\G Vu)^2=(V\G u)^2+(Hu)^2+\la^2 I^2 (Vu)^2\]
\[-2V\G u\cdot Hu+2V\G u\cdot\la IVu-2\la IVu\cdot Hu.\]
Thus:
\[(\G Vu)^2=(V\G u)^2+(Hu)^2+\la^2 I^2 (Vu)^2\]
\[-2V\G u\cdot Hu+2\G Vu\cdot\la IVu-2\la^2 I^2(Vu)^2.\]
Integrating this equation and
\[2\la IVu\cdot \G(V u)=\G((Vu)^{2}\la I)-(Vu)^{2}\cdot \G(\la I)\]
and combining the outcomes with
(\ref{id1}) we arrive at the final integral identity:

\begin{Theorem}
\begin{equation}\label{id}
\int_{SM}(\G Vu)^2\,d\mu-\int_{SM}Q(Vu)^2\,d\mu
=\int_{SM}(V\G u)^2\,d\mu-\int_{SM}(\G u)^2\,d\mu,
\end{equation}
where $Q:={\mathbb K}-\la^2 I^2-\G(\la I)$.
\end{Theorem}

When the Finsler metric is Riemannian (i.e. $I=J=0$), the identity (\ref{id}) is exactly
identity (8) in \cite{DP}.

If $\G u=h(x)+\theta_x(v)$, then one can see that
the right-hand side of (\ref{id}) is nonpositive.  Indeed, since $V\G(u)=V\theta$ we have:
\[\int_{SM}(V\G u)^2\,d\mu-\int_{SM}(\G u)^2\,d\mu
=\int_{SM}(V\theta)^2\,d\mu-\int_{SM}\theta^2\,d\mu-2\int_{SM}h\theta\,d\mu
-\int_{SM}h^{2}\,d\mu.\]
With a bit of work one can see that the linearity of $\theta$ in $v$ implies:
\[\int_{SM}(V\theta)^2\,d\mu=\int_{SM}\theta^2\,d\mu,\]
\[\int_{SM}h\theta\,d\mu=0.\]
This will follow from Lemma \ref{int-nabla}, which holds in any
dimension.

\subsection{Jacobi equation}

For $\zeta\in T(SM)$ write
\[d\phi_{t}(\zeta)=x(t)\G+y(t)H+z(t) V,\]
where $x(t)$, $y(t)$ and $z(t)$ are smooth functions. Equivalently,
\[\zeta=x(t)d\phi_{-t}(\G)+y(t)d\phi_{-t}(H)+z(t)d\phi_{-t}(V).\]
If we differentiate the last equality with respect to $t$ we obtain:
\[0=\dot{x}\G+\dot{y}H+y[\G,H]+\dot{z}V+z[\G,V].\]
Using the bracket relations and regrouping we have:
\[0=(\dot{x}-\la y)\G+(\dot{y}-z-\la Iy)H+(\dot{z}+y{\mathbb K}+z\la I)V,\]
hence
\begin{align*}
&\dot{x}=\la y;\\
&\dot{y}=z+\la Iy;\\
&\dot{z}=-\la Iz-{\mathbb K}y.\\
\end{align*}
From these equations we get:

\begin{equation}
\ddot{y}+Qy=0.
\label{Jac}
\end{equation}

\subsection{Index form}

\begin{Lemma}\label{indexl2d}
If $\phi$ is Anosov, then for every
closed magnetic geodesic $\gamma:[0,T]\to M$ and every
smooth function $z:[0,T]\to \mathbb R$ such that $z(0)=z(T)$
and $\dot z(0)=\dot z(T)$ we have
$$
{\mathbb I}:=\int_0^T \left\{\dot{z}^2
-Qz^2\right\}\,dt\ge0
$$
with equality if and only if $z\equiv 0$.
\end{Lemma}

Using (\ref{Jac}) the proof of this lemma is quite similar to the proof of Lemma 3.3 in \cite{DP}.
The proof of the lemma in any dimension is given in Lemma \ref{crux2}. A key ingredient
is the transversality of the weak stable (or unstable) bundle of $\phi$ with respect
to the vertical distribution, which implies the abscence of conjugate points.

\subsection{End of the proof of Theorem B for surfaces}

Set $\psi:=V(u)$.
The last lemma, applied to the function $z=\psi(\gamma)$, yields
\begin{equation}\label{in_g}
\int_\gamma \left\{(\G\psi)^2
-Q\psi^2\right\}\,dt\ge 0
\end{equation}
for every closed magnetic geodesic $\gamma$.
Since the flow is Anosov, the invariant measures supported on
closed orbits are dense in the space of all invariant measures
on $SM$.
Therefore, the above yields
$$
\int_{SM}\left\{(\G\psi)^2
-Q\psi^2\right\}\,d\mu\ge0.
$$
Combining this with the fact that the right hand side of (\ref{id}) is non-positive, we find that
\begin{equation}\label{in_sm}
\int_{SM}\left\{(\G \psi)^2
-Q\psi^2\right\}\,d\mu=0.
\end{equation}

By the non-negative version of the Liv\v sic theorem, proved independently by M. Pollicott and R. Sharp and by
A. Lopes and P. Thieullen (see \cite{LT, PS}),
we conclude from (\ref{in_g}) and (\ref{in_sm}) that
$$
\int_\gamma\left\{(\G \psi)^2
-Q\psi^2\right\}\,dt=0
$$
for every closed magnetic geodesic $\gamma$.
Applying again Lemma \ref{indexl2d}, we see that
$\psi$ vanishes on all closed magnetic geodesics.
Since the latter are dense in $SM$, the function $\psi$
vanishes on all of $SM$, as required.

\section{Proof of Theorem B}

\subsection{Differential identities of Finsler geometry}\label{basicf}
Henceforth $M$ is a closed $n$-dimensional manifold and
$F$ is a Finsler metric on~$M$.

Let $\pi:TM\setminus\{0\}\to M$ be the natural projection, and let
$\beta^r_s M:=\pi^*\tau^r_s M$ denote the bundle of semibasic
tensors of degree $(r,s)$, where $\tau^r_s M$ is the bundle of
tensors of degree $(r,s)$ over $M$. Sections of the bundles
$\beta^r_sM$ are called semibasic tensor fields and the space of
all smooth sections is denoted by $C^{\infty}(\beta^r_sM)$. For
such a field $T$, the coordinate representation
$$
T=(T^{i_1\dots i_r}_{j_1\dots j_s})(x,y)
$$
holds in the domain of a standard local coordinate system $(x^i,y^i)$
on $TM\setminus\{0\}$
associated with a local coordinate system $(x^i)$ in $M$.
Under a change of a local coordinate system, the components of
a semibasic tensor field are transformed by  the same formula
as those of an ordinary tensor field on $M$.

Every ``ordinary'' tensor field on $M$ defines
a semibasic tensor field by the rule $T\mapsto T\circ\pi$, so that
the space of tensor fields on $M$ can be treated as embedded in the space
of semibasic tensor fields.

Let $(g_{ij})$ be the fundamental tensor,
$$
g_{ij}(x,y)=\frac12[F^2]_{y^iy^j}(x,y),
$$
and let $(g^{ij})$ be the contravariant fundamental tensor,
\begin{equation}\label{g-1}
g_{ik}g^{kj}=\delta_i^j.
\end{equation}

In the usual way, the fundamental tensor
defines the inner product $\langle\cdot,\cdot\rangle$
on $\beta^1_0 M$, and we put $|U|^2=\langle U,U\rangle$.

Let
$$
\mathbf G =y^i\de{}{x^i}-2G^i\de{}{y^i}
$$
be the spray induced by $F$. Here
$G^i$ are the geodesic coefficients \cite[(5.7)]{She},
$$
G^i(x,y)=\frac14g^{il}
\left\{2\de{g_{jl}}{x^k}-\de{g_{jk}}{x^l}\right\}y^jy^k.
$$

Let
\begin{equation*}
T(TM\setminus\{0\})
=\mathcal H TM\oplus \mathcal V TM
\end{equation*}
be the decomposition of $T(TM\setminus\{0\})$
into horizontal and vertical vectors. Here
$$
\mathcal H TM=\operatorname{span}\left\{\del{}{x^i}\right\},
\quad
\mathcal V TM =\operatorname{span}\left\{\de{}{y^i}\right\},
$$
with
$$
\del{}{x^i}=\de{}{x^i}-N^j_i\de{}{y^j}
$$
and
$$
N^i_j=\de{G^i}{y^j}.
$$

Let
$$
\nabla:C^{\infty}(T(TM))\times C^\infty(\pi^*TM)\to
C^{\infty}(\pi^*TM)
$$
be the Chern connection,
$$
\nabla_{\hat X}U=\left\{dU^i(\hat X)+U^j\omega_j^i(\hat X)\right\}\de{}{x^i},
$$
where
$$
\omega^i_j=\Gamma^i_{jk}dx^k
$$
are the connection forms. Recall that
\begin{equation}\label{nij}
N^i_j=\Gamma^i_{jk}y^k.
\end{equation}

Given a function $u\in C^\infty(TM\setminus\{0\})$, we put
$$
u_{|k}:=\del u{x^k},
\quad u_{\cdot k}:=\de u{y^k}
$$
and, given a semibasic vector field $U=(U^i)\in
C^\infty(\beta^1_0M)$, put
$$
U^i_{|k}:=\left(\nabla_{\del{}{x^k}} U\right)^i,
\quad
U^i_{\cdot k}:=\left(\nabla_{\de{}{y^k}} U\right)^i.
$$

We have
$$
u_{|k}=\de u{x^k}-\Gamma^p_{kq}y^q\de u{y^p},
\quad
u_{\cdot k}=\de u{y^k},
$$
and
$$
U^i_{|k}=\de {U^i}{x^k}-\Gamma^p_{kq}y^q\de {U^i}{y^p}
+\Gamma^i_{kp}U^p,
\quad
U^i_{\cdot k}=\de{U^i}{y^k}.
$$

In the usual way, we extend these formulas to higher order
tensors:
\begin{multline*}
T^{i_1\dots i_r}_{j_1\dots j_s|k}
=\de{}{x^k}T^{i_1\dots i_r}_{j_1\dots j_s}
-\Gamma^p_{kq}y^q\de{}{y^p}T^{i_1\dots i_r}_{j_1\dots j_s}
\\
+\sum_{m=1}^r\Gamma^{i_m}_{kp} T^{i_1\dots i_{m-1}pi_{m+1}\dots
i_r}_{j_1\dots j_s} -\sum_{m=1}^s\Gamma^{p}_{kj_m} T^{i_1\dots
i_r}_{j_1\dots j_{m-1}pj_{m+1}\dots j_s}
\end{multline*}
and
$$
T^{i_1\dots i_r}_{j_1\dots j_s\cdot k}
=\de{}{y^k}T^{i_1\dots i_r}_{j_1\dots j_s}.
$$

We define the operators
$$
\nabla_{|}:C^\infty(\beta^r_sM)\to C^\infty(\beta^r_{s+1}M),
\quad
\nabla_{\cdot}:C^\infty(\beta^r_sM)\to C^\infty(\beta^r_{s+1}M)
$$
by
$$
(\nabla_{|} T)^{i_1\dots i_r}_{j_1\dots j_s k}
=\nabla_{|k} T^{i_1\dots i_r}_{j_1\dots j_{s}}
:=T^{i_1\dots i_r}_{j_1\dots j_s|k}
$$
and
$$
(\nabla_{\cdot} T)^{i_1\dots i_r}_{j_1\dots j_{s}k}
=\nabla_{\cdot k} T^{i_1\dots i_r}_{j_1\dots j_{s}}
=T^{i_1\dots i_r}_{j_1\dots j_s\cdot k}.
$$

For convenience, we also define $\nabla^{|}$ and $\nabla^{\cdot}$
by
$$\nabla^{|i}=g^{ij} \nabla_{|j},\quad
\nabla^{\cdot i}=g^{ij}\nabla_{\cdot j}.
$$

In the case of Riemannian manifolds, the above operators
were denoted in \cite{PSh, Sha} by $\overset h \nabla$
and $\overset v \nabla$ respectively.

Given a function $u\in C^\infty(TM\setminus\{0\})$, note that
$\nabla^{\cdot}u=0$ if and only if $u$ does not depend on $y$.

Equivalently, the above can be described as follows.
In a natural way, the connection $\nabla$ on $\beta^1_0 M=\pi^*TM$ defines
a connection on the dual bundle $\beta^0_1=\pi^*T^*M$,
as well as connections on the tensor product bundles $\beta^r_sM$
for all $r$ and $s$.
Then for $T\in C^\infty(\beta^r_s M)$ we have
$$
\nabla_{|k} T=\nabla_{\del{}{x^k}} T,
\quad \nabla_{\cdot k} T=\nabla_{\de{}{y^k}} T.
$$

This shows also that $\nabla_{|}$ and $\nabla_{\cdot}$
are compatible with tensor products and contractions.

Note that
$$
g_{ij\cdot k}=2C_{ijk},
\quad
g^{ij}_{\cdot k}=-2g^{il}g^{jm}C_{lmk},
$$
where
$$
C_{ijk}=\frac14[F^2]_{y^iy^jy^k}
$$
is the Cartan tensor of $F$.

Also, note that
the fundamental tensor is parallel with respect to $\nabla_{|}$:
\begin{equation}\label{nablahg}
g_{ij|k}=0\quad g^{ij}_{|k}=0.
\end{equation}
Indeed, using (5.29) of \cite{She}, we see that
$$
g_{ij|k}=\de{g_{ij}}{x^k}
-\Gamma^p_{kq}y^q\de{g_{ij}}{y^p}
-\Gamma^p_{ki}g_{pj}-\Gamma^p_{kj}g_{ip}
=2C_{ipj}N^p_k-2\Gamma^p_{kq}y^qC_{ijp}=0,
$$
while the second identity is obtained by differentiating (\ref{g-1}).

By \cite[Lemma 5.2.1]{She}
$$
F_{|k}=0.
$$

On the other hand, for $(x,y)\in SM$
\begin{equation}\label{fk}
F_{\cdot k}=y_k=g_{kj}y^j.
\end{equation}
Indeed, using homogeneity we have
$$
FF_{y^k}=\frac12[F^2]_{y^k}=\frac12[F^2]_{y^k y^j}y^j
=g_{kj} y^j.
$$
However, $F=1$ on $SM$, which gives \eqref{fk}.

A straightforward computation shows also that
$$
y^i_{|k}=0,\quad y^i_{\cdot k}=\delta^i_j.
$$

Let $P$ denote the Chern curvature tensor and $R$ denote
the Riemann curvature tensor (see \cite[(8.12), (8.13)]{She}):
$$
P^i_{jkl}=-\de{\Gamma^i_{jk}}{y^l},
$$
$$
R^i_{jkl}=\de{\Gamma^i_{jl}}{x^k}
-\de{\Gamma^i_{jk}}{x^l}
+\de{\Gamma^i_{jk}}{y^m}N^m_l
-\de{\Gamma^i_{jl}}{y^m}N^m_k
+\Gamma^m_{jl}\Gamma^i_{mk}
-\Gamma^m_{jk}\Gamma^i_{ml},
$$
and put (see \cite[p. 127]{She})
$$
P^i_{kl}=y^j P^i_{jkl},
$$
$$
R^i_{kl}=y^j R^i_{jkl}.
$$

Note that
$$
R^i_k=R^i_{kl}y^l
$$
corresponds to the Riemann curvature operator
$$
\mathbf R_y(V)=(R^i_{k}V^k),
$$
while
\begin{equation}\label{yp}
y^kP^i_{kl}=0.
\end{equation}

\begin{Lemma}
If $u\in C^\infty(TM\setminus\{0\})$, then
\begin{align}\label{vv}
u_{\cdot l\cdot k}- u_{\cdot k\cdot l}&=0,
\\
\label{hv}
u_{|l\cdot k}- u_{\cdot k|l}
&=P^i_{lk}u_{\cdot i},
\\
\label{hh}
u_{|l|k}- u_{|k|l}
&=R^i_{lk}u_{\cdot i}.
\end{align}
\end{Lemma}

\begin{proof}
(\ref{vv}) is trivial.

Next,
$$
u_{|l\cdot k}=\de{}{y^k}\left(\de u{x^l}
-\Gamma^i_{lj}y^j\de u{y^i}\right)
=\frac{\partial^2 u}{\partial y^k\partial x^l}
-\de{\Gamma^i_{lj}}{y^k}y^j\de u{y^i}
-\Gamma^i_{lk}\de u{y^i}
-\Gamma^i_{lj}y^j\frac{\partial^2 u}{\partial y^k\partial y^i},
$$
whereas
$$
u_{\cdot k|l}
=\left(\de{}{x^l}
-\Gamma^i_{lj}y^j\de{}{y^i}\right) u_{\cdot k}
-\Gamma^i_{lk}u_{\cdot i}
=
\frac{\partial^2 u}{\partial x^l\partial y^k}
-\Gamma^i_{lj}y^j\frac{\partial^2 u}{\partial y^i\partial y^k}
-\Gamma^i_{lk}\de u{y^i}.
$$
Taking the difference, we come to (\ref{hv}).

Further,
\begin{multline*}
u_{|l|k}
=\left(\de{}{x^k}-\Gamma^m_{ks}y^s\de{}{y^m}\right)
u_{|l}-\Gamma^m_{kl}u_{|m}
\\
=\left(\de{}{x^k}-\Gamma^m_{ks}y^s\de{}{y^m}\right)
\left(\de{u}{x^l}-\Gamma^i_{lj}y^j\de u{y^i}\right)
-\Gamma^m_{kl}\left(\de{u}{x^m}-\Gamma^i_{mj}y^j\de u{y^i}\right)
\\
=\frac{\partial^2 u}{\partial x^k\partial x^l}
-\Gamma^m_{ks}y^s\frac{\partial^2 u}{\partial y^m\partial x^l}
-\de{\Gamma^i_{lj}}{x^k}y^j\de u{y^i}
-\Gamma^i_{lj}y^j\frac{\partial^2u}{\partial x^k\partial y^i}
\\
+\Gamma^m_{ks}y^s\de{\Gamma^i_{lj}}{y^m}y^j\de u{y^i}
+\Gamma^m_{ks}y^s\Gamma^i_{lm}\de u{y^i}
+\Gamma^m_{ks}y^s\Gamma^i_{lj}y^j\frac{\partial^2 u}{\partial y^m\partial y^i}
-\Gamma^m_{kl}\de u{x^m}
+\Gamma^m_{kl}\Gamma^i_{mj}y^j\de u{y^i}.
\end{multline*}

Using (\ref{nij}), rearranging, and appropriately renaming
indices, we obtain
\begin{multline*}
u_{|l|k}
=\frac{\partial^2 u}{\partial x^k\partial x^l}
-N^m_k\frac{\partial^2 u}{\partial y^m\partial x^l}
-N^i_l\frac{\partial^2u}{\partial x^k\partial y^i}
+N^m_k N^i_l\frac{\partial^2 u}{\partial y^m\partial y^i}
-\Gamma^m_{kl}\de u{x^m}
\\
-\left(\de{\Gamma^i_{lj}}{x^k}
-\de{\Gamma^i_{lj}}{y^m}N^m_k
-\Gamma^m_{kj}\Gamma^i_{lm}
-\Gamma^m_{kl}\Gamma^i_{mj}\right)y^j\de u{y^i}.
\end{multline*}

Alternating with respect to $k$ and $l$, we
come to (\ref{hh}).
\end{proof}

\subsection{Integral identities of Finsler geometry}
We will derive the Gauss--Ostrogradski\u\i{} formulas
for vertical and horizontal divergences like those for Riemannian manifolds
in \cite[Section 3.6]{Sha}. We proceed along the lines of \cite{Sha}.

Given a vector field $U=(U^i)\in C^\infty(\beta^1_0M)$,
the vertical divergence and the horizontal divergence are defined by
$$
\divv U=U^i_{\cdot i},
\quad
\divh U=U^i_{|i}.
$$

Let
$$
\mathbf I(U)=g^{ij}C_{ijk}U^k
$$
be the mean Cartan torsion \cite[p. 108]{She}, and let
$$
\mathbf J(U)=g^{ij}L_{kij}U^k
$$
be the mean Landsberg curvature
\cite[p. 116]{She}).
Here $L$ is the Landsberg tensor,
related to the Chern curvature tensor
as follows \cite[(8.27)]{She}:
\begin{equation}\label{l-p}
L_{ijk}=-g_{im}P^m_{jk}.
\end{equation}

Let
$$
dV^{2n}=\det (g_{ij})\,dx^1\dots dx^ndy^1\dots dy^n
$$
be the Liouville volume form on $TM\setminus\{0\}$.

Consider the following set of local forms on $TM\setminus\{0\}$
$$
\ov_k=(-1)^{n+k-1}g\,dx\wedge dy^1\wedge
\dots\wedge \widehat{dy^k}\wedge \dots \wedge  dy^n,
$$
\begin{multline*}
\oh_k=g\big[(-1)^{k-1}\,dx^1\wedge \dots\wedge \widehat{dx^k}\wedge
\dots \wedge dx^n \wedge dy
\\
+\sum_{j=1}^n (-1)^{n+j}\Gamma^j_{kl} y^l\,dx\wedge dy^1\wedge
\dots\wedge \widehat{dy^j}\wedge \dots\wedge  dy^n\big],
\end{multline*}
where $g=\det(g_{ij})$, $dx=dx^1\wedge \dots\wedge  dx^n$,
$dy=dy^1\wedge \dots\wedge  dy^n$, and the symbol $\;$$\widehat{}$$\;$ over a factor
means that this factor is omitted.

\begin{Lemma}
Given a semibasic vector field $U=(U^k)$, the set of local forms
$U^k\ov_k$ defines a global differential form on $TM\setminus\{0\}$.
Similarly, the set of local forms
$U^k\oh_k$ defines a global differential form on $TM\setminus\{0\}$.
Moreover,
\begin{equation}\label{dov}
d(U^k\ov_k)=(\divv U+2\mathbf I(U))\,dV^{2n},
\end{equation}
\begin{equation}\label{doh}
d(U^k\oh_k)=(\divh U-\mathbf J(U))\,dV^{2n}.
\end{equation}
\end{Lemma}

\begin{proof}
$$
d\ov_k=\de g{y^k}\,dx\wedge dy=g^{ij}\de{g_{ij}}{y^k}g\,dx\wedge dy
=2g^{ij}C_{ijk}\,dV^{2n}.
$$
Therefore,
$$
d(U^k\ov_k)=\de{U^k}{y^k}g\,dx\wedge dy
+2U^k g^{ij}C_{ijk}\,dV^{2n},
$$
which coincides with (\ref{dov}).

Next,
$$
\begin{aligned}
d\oh_k
&=\de g{x^k}\,dx\wedge dy
-\left(\de{g}{y^j}\Gamma^j_{kl}y^l
+g\de{\Gamma^j_{kl}}{y^j}y^l
+g\Gamma^j_{kj}\right)\,dx\wedge dy
\\
&=\left(g^{ij}\de{g_{ij}}{x^k}
-g^{km}\de{g_{km}}{y^j}\Gamma^j_{kl}y^l
-\de{\Gamma^j_{kl}}{y^j}y^l
-\Gamma^j_{kj}\right)g\,dx\wedge dy
\\
&=\left(g^{ij}\de{g_{ij}}{x^k}
-2g^{km}C_{kmj}N^j_k
-\Gamma^j_{kj}
+P^j_{kj}\right)\,dV^{2n}
=(\Gamma^j_{kj}+P^j_{kj})\,dV^{2n}.
\end{aligned}
$$
Here we have used the equality \cite[(5.29)]{She}
$$
\de{g_{jl}}{x^m}=g_{kl}\Gamma^k_{jm}+g_{kj}\Gamma^k_{lm}+2C_{jkl}N^k_m.
$$

Consequently,
$$
\begin{aligned}
d(U^k\oh_k)&=\de{U^k}{x^k}g\,dx\wedge dy
-\de{U^k}{y^j}g\Gamma^j_{kl}y^l\,dx\wedge dy
+U^k(\Gamma^j_{kj}+P^j_{kj})\,dV^{2n}
\\
&=\left\{\left(\de{U^k}{x^k}-\Gamma^j_{kl}y^l\de{U^k}{y^j}
+\Gamma^j_{kj}U^k\right)+P^j_{kj} U^k\right\}\,dV^{2n},
\end{aligned}
$$
which coincides with (\ref{doh}) in view of (\ref{l-p})
and the symmetry of the Landsberg tensor.
\end{proof}

Let $SM=\{(x,y)\in TM\mid F(y)=1\}$
be the unit sphere bundle. The restriction of the form
$y^k\ov_k$ to $SM$ gives rise to
the Liouville measure $d\mu$ of $SM$.

\begin{Theorem}\label{g-o}
Let $U\in C^\infty(\beta^1_0M)$ be a semibasic vector field
positively homogeneous of degree $\lambda$ in $y$.
Then the following Gauss--Ostrogradski\u\i{} formulas hold:
\begin{align}\label{divv}
\int_{SM}\divv U\,d\mu&=\int_{SM}((\lambda+n-1)\langle U,y\rangle
-2\mathbf I(U))\,d\mu,
\\
\label{divh}
\int_{SM}\divh U\,d\mu&=\int_{SM}\mathbf J(U)\,d\mu.
\end{align}
\end{Theorem}

These formulas follow easily from (\ref{dov})--(\ref{doh})
by integration.

\begin{Lemma}\label{int-nabla}
\begin{enumerate}
\item Let $\psi\in C^{\infty}(TM)$ be a function which depends linearly on $y$. Then
\[\int_{SM}\psi\,d\mu=0.\]
\item Let $\phi\in C^\infty(TM\setminus\{0\})$ be such that
$\phi=\varphi_0 F+\psi$, where $\varphi_0$
is independent of $y$ while $\psi$ depends linearly on $y$. Then
$$
\int_{SM}|\nabla^{\cdot}\phi|^2\,d\mu
=\int_{SM}(\varphi_0^2+n\psi^2)\,d\mu.
$$
\end{enumerate}
\end{Lemma}

\begin{proof} To prove (1) let $\psi=\Psi_k y^k$,
where $\Psi$ is a covector field on $M$.
Put $U^i=g^{ij}\Psi_j$ and apply (\ref{divv}) to get
$$
(n-1)\int_{SM}\psi\,d\mu
=(n-1)\int_{SM}\langle U,y\rangle \,d\mu
=\int_{SM}(\divv U+2\mathbf I(U))\,d\mu.
$$
Now,
$$
\divv U=(g^{ij}\Psi_j)_{\cdot i}
=g^{ij}_{\cdot i}\Psi_j+g^{ij}\Psi_{j\cdot i}
=-2g^{il}g^{jm}C_{lmi}\Psi_j=-2\mathbf I(U),
$$
which implies (1).

To prove (2) note that since
$\nabla_{\cdot}\phi=\varphi_0\nabla_{\cdot}F+\nabla_{\cdot}\psi$, we have
\begin{equation*}
|\nabla_{\cdot}\phi|^2
=\varphi_0^2|\nabla^{\cdot}F|^2
+2\varphi_0\langle\nabla_{\cdot}F,\nabla_{\cdot}\psi\rangle
+|\nabla^{\cdot}\psi|^2.
\end{equation*}

Next,
\begin{equation*}
\begin{aligned}
|\nabla^{\cdot}\psi|^2
&=g^{ij}\psi_{\cdot i}\psi_{\cdot j}
=(\psi g^{ij}\psi_{\cdot i})_{\cdot j}
-\psi g^{ij}_{\cdot j}\psi_{\cdot i}
-\psi g^{ij}\psi_{\cdot i\cdot j}
\\
&=\divv(\psi\nabla^{\cdot}\psi)
+2\mathbf I(\psi\nabla^{\cdot}\psi),
\end{aligned}
\end{equation*}
because $\psi_{\cdot i\cdot j}=0$.

Thus, on $SM$ we get
$$
|\nabla^{\cdot}\phi|^2
=\varphi_0^2+2\varphi_0\psi
+\divv(\psi\nabla^{\cdot}\psi)
+2\mathbf I(\psi\nabla^{\cdot}\psi).
$$

Integrating  and using (\ref{divv}), we receive
$$
\int_{SM}|\nabla^{\cdot}\phi|^2\,d\mu
=\int_{SM}(\varphi_0^2+2\varphi_0\psi+n\psi\langle
\nabla^{\cdot}\psi,y\rangle)\,d\mu =\int_{SM}
(\varphi_0^2+2\varphi_0\psi+n\psi^2)\,d\mu.
$$
Since by (1)
\[\int_{SM}\varphi_0\psi\,d\mu=0\]
the proof of (2) is complete.

\end{proof}

\subsection{Identities for the magnetic flow}\label{idloc}
Let
$$
\{dx^i,\ \delta y^j= d y^j+N^j_k\, dx^k\}
$$
be a local basis for $T^*(TM\setminus\{0\})$  dual
to the local basis
$\left\{\del{}{x^i},\ \de{}{y^j}\right\}$ for $T(TM\setminus\{0\})$.
The Legendre transform $\ell_{F}:TM\setminus\{0\}\to T^*M\setminus\{0\}$
associated with the Lagrangian $\frac{1}{2}F^2$ is a diffeomorphism and
$\omega_{0}:=\ell_{F}^{*}(-d\la)$ defines a symplectic form on $TM\setminus\{0\}$,
where $\la$ is the Liouville 1-form on $T^*M$.
In local coordinates $(x,y)$, $\ell_{F}$ is simply the map
\[(y^{j})\mapsto (y_{j}).\]
The canonical 1-form is $\la=y_{i}dx^{i}$ and
$\ell_{F}^{*}\la=g_{ij}y^{j}dx^{i}$. From this, a calculation
shows that
$$
\omega_0=g_{ij}dx^i\wedge\delta y^j.
$$

Let $H:TM\setminus \{0\}\to \mathbb R$ be defined by
$$
H=\frac 12 F^2.
$$
The Hamiltonian flow of $H$ with respect to $\omega_0$
gives rise to the geodesic flow of the Finsler manifold $(M,F)$.

Let $\Omega$ be a closed 2-form on $M$ and consider
the new symplectic form $\omega$ defined as
$$
\omega_0+\pi^*\Omega.
$$
The Hamiltonian flow of $H$ with respect to
$\omega_{0}+\pi^{*}\Omega$ gives rise to a flow
$\phi_t:TM\setminus\{0\}\to TM\setminus\{0\}$, called {\em
magnetic flow} or {\em twisted geodesic flow}.

The form $\Omega$, regarded as an antisymmetric tensor field
$(\Omega_{ij})\in C^\infty(\tau^0_2 M)$, gives rise to a corresponding
semibasic tensor field. We
define the {\it Lorentz force} $Y\in C^\infty(\beta^1_1M)$ by
\begin{equation}\label{lorentz}
Y^i_j(x,y)=\Omega_{jk}(x)g^{ik}(x,y).
\end{equation}

We also define
$$
Y(U)=(Y^i_j U^j).
$$

Note that $Y$ is skew symmetric with respect to $g$:
$$
\langle Y(U),V\rangle=-\langle U,Y(V)\rangle.
$$

Let $\mathbf G_M$ be the generator of the magnetic flow.
Straightforward calculations show that
\begin{equation}\label{spray}
\mathbf G_M(x,y)=y^i\del{}{x^i}+y^iY^j_i\de{}{y^j}.
\end{equation}

It is easily seen that every integral curve of
$\mathbf G_M$ is a curve of the form
$t\mapsto \dot\gamma(t)\in TM$
which satisfies the equation
$$
D_{\dot \gamma}\dot\gamma=Y_{\dot\gamma(t)}(\dot\gamma),
$$
where the covariant derivative $D$ is the one determined by the
Chern connection. Alternatively we could write:
$$
\ddot\gamma^i(t)
+\Gamma^i_{jk}(\dot\gamma(t))\dot\gamma^j(t)\dot\gamma^k(t)
=Y^i_j(\dot\gamma(t))\dot\gamma^j(t).
$$
A curve $\gamma$, satisfying this equation, is referred to as a
magnetic geodesic.

If $u\in C^\infty(TM\setminus\{0\})$, then by (\ref{spray})
$$
\mathbf G_M u(x,y)=y^i\left(\del u{x^i}+Y^j_i\de u{y^j}\right)
=y^i(u_{|i}+Y^j_iu_{\cdot j}).
$$

Since the Hamiltonian flow $\phi_t$ preserves the level sets of $H$,
the magnetic flow preserves $SM$
and the vector field $\mathbf G_M$ is tangent to $SM$.

Suppose that for a smooth function $u:SM\to \mathbb R$
we have
$$
\mathbf G_M u=\varphi.
$$
Extend $u$ to a positively homogeneous function (of degree $0$)
on $TM\setminus\{0\}$, denoting the extension by $u$ again.

For $(x,y)\in TM$, define
$$
\mathbf X u=y^i(u_{|i}+FY^j_iu_{\cdot j}).
$$

Then on $TM\setminus\{0\}$
we have
$$
\mathbf Xu=\phi,
$$
where $\phi$ is the positively homogeneous extension of $\varphi$ to
$TM\setminus\{0\}$ of degree $1$.

Given $T=(T^{i_1\dots i_r}_{j_1\dots j_s})\in C^\infty(\beta^r_s M)$,
put
$$
T^{i_1\dots i_r}_{j_1\dots j_s: k}
=T^{i_1\dots i_r}_{j_1\dots j_s|k}
+FY^j_k T^{i_1\dots i_r}_{j_1\dots j_s\cdot j}.
$$

Straightforward calculations show that for $(x,y)\in SM$
\begin{align}
g_{ij:k}&=2Y^s_k C_{ijs}, \label{gijk}
\\
g^{ij}_{:k}&=-2Y^s_k g^{il}g^{jm}C_{lms}, \nonumber
\\
y^i_{:k}&=Y^i_k. \nonumber
\end{align}

It is also useful to note that differentiating (\ref{lorentz}) yields
\begin{equation}\label{yijk}
Y^i_{j\cdot k}=-2Y^m_j g^{il}C_{lmk}=g^{is}_{:j}g_{sk}.
\end{equation}

\begin{Lemma}\label{com}
If $u\in C^\infty(TM\setminus\{0\})$, then
for $(x,y)\in SM$ we have
\begin{align}\label{mv}
u_{:l\cdot k}-u_{\cdot k:l}
&=\tilde P^i_{lk}u_{\cdot i},
\\
\label{mm}
u_{:l:k}-u_{:k:l}
&=\tilde R^i_{lk} u_{\cdot i},
\end{align}
with
\begin{align*}
\tilde P^i_{lk}&=P^i_{lk}
+Y^i_l y_k
+Y^i_{l\cdot k},
\\
\tilde R^i_{lk}&=R^i_{lk}+(Y^i_{l|k}-Y^i_{k|l})
-(P^i_{lm}Y^m_k-P^i_{km}Y^m_l)
\\
&\quad+(Y^j_lY^i_{k\cdot j}-Y^j_kY^i_{l\cdot j})
+y_s(Y^s_k Y^i_l-Y^s_l Y^i_k).
\end{align*}
\end{Lemma}

\begin{proof}
We have
$$
u_{:l\cdot k}=(u_{|l}+FY^i_lu_{\cdot i})_{\cdot k}
=u_{|l\cdot k}+F_{\cdot k} Y^i_l u_{\cdot i}
+FY^i_{l\cdot k}u_{\cdot i}+FY^i_l u_{\cdot i\cdot k}
$$
whereas
$$
u_{\cdot k:l}
=u_{\cdot k|l}+FY^i_lu_{\cdot k\cdot i}.
$$
Thus, for $(x,y)\in SM$
$$
u_{:l\cdot k}-u_{\cdot k:l}
=(u_{|l\cdot k}-u_{\cdot k|l})
+y_kY^i_lu_{\cdot i}+Y^i_{l\cdot k}u_{\cdot i}.
$$

Using (\ref{hv}), we come to (\ref{mv}).

Further,
\begin{multline*}
u_{:l:k}=u_{:l|k}+FY^j_k u_{:l\cdot j}
=(u_{|l}+FY^j_lu_{\cdot j})_{|k}
+FY^j_k(u_{|l}+FY^s_l u_{\cdot s})_{\cdot j}
\\=u_{|l|k}+FY^j_{l|k}u_{\cdot j}
+FY^j_l u_{\cdot j|k}
+FY^j_ku_{|l\cdot j}
\\
+FY^j_kF_{\cdot j} Y^s_lu_{\cdot s}
+F^2Y^j_kY^s_{l\cdot j} u_{\cdot s}
+F^2Y^j_kY^s_l u_{\cdot s\cdot j}.
\end{multline*}

Thus, for $(x,y)\in SM$
\begin{multline*}
u_{:l:k}-u_{:k:l}
=(u_{|l|k}-u_{|k|l})
+(Y^j_{l|k}-Y^j_{k|l})u_{\cdot j}
\\
+Y^j_l (u_{\cdot j|k}-u_{|k\cdot j})
+Y^j_k(u_{|l\cdot j}-u_{\cdot j|l})
+(Y^j_k Y^s_l-Y^j_l Y^s_k)y_ju_{\cdot s}
\\
+(Y^j_kY^s_{l\cdot j}-Y^j_lY^s_{k\cdot j})u_{\cdot s}
+(Y^j_kY^s_l-Y^j_lY^s_k) u_{\cdot s\cdot j}.
\end{multline*}
Using (\ref{hh}), (\ref{hv}) and renaming indices, we come to (\ref{mm}).
\end{proof}

Given $U\in C^\infty(\beta^1_0M)$ and $u\in
C^\infty(TM\setminus\{0\})$, define
$$
\divm U=U^i_{:i},\quad
\nabla^{:}u=(u^{:i})=(g^{ij} u_{:j}).
$$

\begin{Lemma}\label{pes-n}
The following holds on $SM$ {\rm (The Pestov identity):}
\begin{align}
\nonumber 2\langle\nabla^{:} u,\nabla^{\cdot}(\mathbf X u)\rangle
&=|\nabla^{:} u|^2 +\mathbf X (\langle\nabla^{:} u,\nabla^{\cdot}
u\rangle) -\divm((\mathbf X u)\nabla^{\cdot} u) +\divv((\mathbf X
u)\nabla^{:} u)
\\
\nonumber
&\quad-\langle \tilde {\mathbf R}_y(\nabla^{\cdot} u),\nabla^{\cdot} u\rangle
+\langle Y(\nabla^{\cdot}u),\nabla^{:}u\rangle
\\
\label{pestov}
&\quad+2\mathbf I((\mathbf X u)\nabla^{:} u)
+\mathbf J((\mathbf X u)\nabla^{\cdot} u).
\end{align}
\end{Lemma}

\begin{proof}
With the above notations, we can write
$$
\mathbf X u= y^i u_{:i}.
$$
Therefore,
\begin{multline}\label{start}
2\langle\nabla_{\cdot}(\mathbf X u),\nabla_{:} u\rangle
-\divv((\mathbf X u)\nabla^{:} u)
=2g^{ij}(\mathbf X u)_{\cdot i} u_{:j}
-((\mathbf X u) g^{ij} u_{:j})_{\cdot_i}
\\
=g^{ij}(\mathbf X u)_{\cdot_i} u_{:j}
-(\mathbf X u) g^{ij}_{\cdot_i} u_{:j}
-(\mathbf X u) g^{ij}u_{:j\cdot i}=I-II-III.
\end{multline}

We rewrite the first term on the right-hand side of (\ref{start}) as follows:
\begin{align*}
I&=g^{ij}(y^k u_{:k})_{\cdot_i} u_{:j}
=g^{ij}(u_{:i}+y^k u_{:k\cdot i}) u_{:j}
\\
&=g^{ij} u_{:i} u_{:j}
+g^{ij}y^k(u_{\cdot i:k}+(u_{:k\cdot i}-u_{\cdot i:k})) u_{:j}
\\
&=|\nabla^{:} u|^2
+y^k(g^{ij} u_{\cdot i} u_{:_j})_{:k}
-y^kg^{ij}_{:k}u_{\cdot i}u_{:j}
-y^k g^{ij}u_{\cdot i} u_{:j:k}
+g^{ij}y^k \tilde P^m_{ki} u_{\cdot m} u_{:j}.
\end{align*}

Note that
$$
y^k(g^{ij} u_{\cdot i} u_{:_j})_{:k}
=\mathbf X (\langle\nabla^{\cdot} u,\nabla^{:} u\rangle),
$$
that
\begin{align*}
g^{ij}y^k \tilde P^m_{ki} u_{\cdot m} u_{:j}
&= g^{ij}y^k(P^m_{ki}+Y^m_k y_i+Y^m_{k\cdot i})u_{\cdot m} u_{:j}
\\
&=\langle Y(y),\nabla^{\cdot}u\rangle\mathbf X u
+y^kg^{mj}_{:k}u_{\cdot m}u_{:j}
\end{align*}
where we have used (\ref{yp}) and (\ref{yijk}), and that
\begin{align*}
y^k g^{ij}u_{\cdot i} u_{:j:k}
&=y^k g^{ij} u_{\cdot i}(u_{:k:j}+(u_{:j:k}-u_{:k:j}))
\\
&=g^{ij} u_{\cdot _i}(y^k u_{:_k})_{:j}
-g^{ij}u_{\cdot i}y^k_{:j}u_{:k}
+y^kg^{ij} u_{\cdot i} \tilde R^m_{jk} u_{\cdot m}
\\
&=\langle\nabla^{\cdot} u,\nabla^{:}(\mathbf X u)\rangle
-\langle Y(\nabla^{\cdot}u),\nabla^{:}u\rangle
+\langle \tilde {\mathbf R}_y(\nabla^{\cdot} u),\nabla^{\cdot} u\rangle.
\end{align*}

Thus,
\begin{multline}\label{start1}
I=|\nabla^{:} u|^2
+\mathbf X (\langle\nabla^{\cdot} u,\nabla^{:} u\rangle)
+\langle Y(\nabla^{\cdot}u),\nabla^{:}u\rangle
-\langle \tilde {\mathbf R}_y(\nabla^{\cdot} u),\nabla^{\cdot} u\rangle.
\\
+\langle Y(y),\nabla^{\cdot}u\rangle\mathbf X u
-\langle\nabla^{\cdot} u,\nabla^{:}(\mathbf X u)\rangle.
\end{multline}

We rewrite the second term on the right-hand side of (\ref{start}) as
\begin{equation}\label{start2}
II=(\mathbf X u) g^{ij}_{\cdot i} u_{:j}
=-2(\mathbf X u) g^{il}g^{jm}C_{lmi} u_{:j}
=-2\mathbf I((\mathbf X u)\nabla^{:} u).
\end{equation}

Finally, we rewrite the third term in (\ref{start}) as
\begin{align*}
III&=(\mathbf X u) g^{ij}u_{:j\cdot i}
=(\mathbf X u) g^{ij}(u_{\cdot i:j} +(u_{:j\cdot i}-u_{\cdot i:j}))
\\
&=((\mathbf X u) g^{ij} u_{\cdot i})_{:j}
-(\mathbf X u)_{:j} g^{ij}u_{\cdot i}
-(\mathbf X u)g^{ij}_{:j} u_{\cdot i}
+(\mathbf X u) g^{ij} \tilde P^m_{ji} u_{\cdot m}.
\end{align*}

Note that
$$
(\mathbf X u) g^{ij} u_{\cdot i})_{:j}
=\divm((\mathbf X u)\nabla^{\cdot} u),
$$
that
$$
(\mathbf X u)_{:j} g^{ij}u_{\cdot i}
=\langle \nabla^{\cdot}u, \nabla^{:}(\mathbf X u)\rangle,
$$
and that
\begin{align*}
(\mathbf X u) g^{ij} \tilde P^m_{ji} u_{\cdot m}
&= (\mathbf X u) g^{ij}(P^m_{ji}+Y^m_j y_i+Y^m_{j\cdot i})u_{\cdot m}
\\
&=-\mathbf J((\mathbf X u)\nabla^{\cdot} u)
+\langle Y(y),\nabla^{\cdot}u\rangle\mathbf X u
+(\mathbf X u)g^{mj}_{:j}u_{\cdot m}
\end{align*}
in view of (\ref{yijk}).

Thus,
\begin{equation}\label{start3}
III=\divm((\mathbf X u)\nabla^{\cdot} u)
-\mathbf J((\mathbf X u)\nabla^{\cdot} u)
+\langle Y(y),\nabla^{\cdot}u\rangle \mathbf X u
-\langle \nabla^{\cdot}u, \nabla^{:}(\mathbf X u)\rangle.
\end{equation}

Inserting (\ref{start1})--(\ref{start3}) in (\ref{start}),
we come to (\ref{pestov}).
\end{proof}

Given a semibasic vector field $V$, define
a new semibasic vector field $\mathbf X V$ by
$$
\mathbf X V^i
=y^kV^i_{:k}.
$$

It easy to see that if $(x,y)\in SM$ and $\gamma$ is a magnetic
geodesic with $\gamma(0)=x$, $\dot\gamma(0)=y$, then
$$
\mathbf X V(x,y)=D_{\dot\gamma} (V\circ \dot\gamma)|_{t=0},
$$
the covariant derivative of the field $V\circ\dot\gamma$ along $\gamma$.

\begin{Lemma}
If $u\in C^\infty(TM\setminus\{0\})$ is positively homogeneous, then
\begin{equation}\label{xnabla}
|\mathbf X(\nabla^{\cdot}u)|^2
=|\nabla^{\cdot} \mathbf X  u|^2
+|\nabla^{:}u|^2
-2\langle\nabla^{:}u,\nabla^{\cdot}(\mathbf X u)\rangle
+\langle Y(y),\nabla^{\cdot}u\rangle^2.
\end{equation}
\end{Lemma}

\begin{proof}
We have
$$
\begin{aligned}
\mathbf X (u^{\cdot i})
&=y^k(g^{ij}u_{\cdot j})_{:k}
=y^kg^{ij}_{:k}u_{\cdot j}
+y^k g^{ij}(u_{:k\cdot j}
-(u_{:k\cdot j}-u_{\cdot j:k}))
\\
&=y^kg^{ij}_{:k}u_{\cdot j}
+g^{ij}(y^k u_{:k})_{\cdot j}
-g^{ij}u_{:j}
-g^{ij}y^k\tilde P^m_{kj} u_{\cdot m}.
\end{aligned}
$$

By (\ref{yp}) and (\ref{yijk})
\begin{align*}
g^{ij}y^k\tilde P^m_{kj} u_{\cdot m}
&=g^{ij}y^k(P^m_{kj} +Y^m_k y_j +Y^m_{k\cdot j})u_{\cdot m}
\\
&=\langle Y(y),\nabla^{\cdot}u\rangle y^i
+y^kg^{mi}_{:k}u_{\cdot m}.
\end{align*}

Thus
\begin{equation*}
\mathbf X(\nabla^{\cdot}u)
=\nabla^{\cdot}(\mathbf Xu)
-\nabla^{:}u
-\langle Y(y),\nabla^{\cdot}u\rangle y.
\end{equation*}

Squaring, we receive
\begin{multline*}
|\mathbf X (\nabla^{\cdot} u)|^2
=|\nabla^{\cdot} \mathbf X  u|^2
+|\nabla^{:}u|^2
+\langle Y(y),\nabla^{\cdot}u\rangle^2
\\
-2\langle\nabla^{\cdot}(\mathbf X u),\nabla^{:}u\rangle
-2\langle Y(y),\nabla^{\cdot}u\rangle
\langle \nabla^{\cdot}(\mathbf Xu),y\rangle
+2\langle Y(y),\nabla^{\cdot}u\rangle\langle\nabla^{:}u,y\rangle
\\
=|\nabla^{\cdot} \mathbf X  u|^2
+|\nabla^{:}u|^2
+\langle Y(y),\nabla^{\cdot}u\rangle^2
-2\langle\nabla^{:}u,\nabla^{\cdot}(\mathbf X u)\rangle
\\
-2\langle Y(y),\nabla^{\cdot}u\rangle\mathbf Xu
+2\langle Y(y),\nabla^{\cdot}u\rangle\mathbf Xu,
\end{multline*}
coming to the sought identity.
\end{proof}

Suppose that we have a kinetic equation on $SM$
$$
\mathbf G_Mu=\varphi.
$$
Extending $u$ to a positively
homogeneous function on $TM\setminus\{0\}$,
the extension denoted by $u$ again, we have on $TM\setminus\{0\}$
$$
\mathbf X u=\phi,
$$
where $\phi$ is the positively homogeneous extension of $\varphi$
of degree $1$.

Combining (\ref{pestov}) and (\ref{xnabla}), we get
\begin{multline*}
|\mathbf X (\nabla^{\cdot} u)|^2
+\mathbf X (\langle\nabla^{:} u,\nabla^{\cdot} u)\rangle
-\divm((\mathbf X  u)\nabla^{\cdot} u)
\\
-\langle \tilde {\mathbf R}_y(\nabla^{\cdot} u),\nabla^{\cdot} u\rangle
+\langle Y(\nabla^{\cdot}u),\nabla^{:} u\rangle
-\langle Y(y),\nabla^{\cdot}u\rangle^2
\\
+2\mathbf I((\mathbf X u)\nabla^{:} u)
+\mathbf J((\mathbf X u)\nabla^{\cdot} u)
\\
=|\nabla^{\cdot} (\mathbf X  u)|^2
-\divv((\mathbf X u)\nabla^{:} u).
\end{multline*}

We integrate this identity over $SM$ against the Liouville measure,
using the flow invariance of the measure and (\ref{divv}):
\begin{multline}\label{integral}
\int_{SM}|\mathbf X (\nabla^{\cdot} u)|^2\,d\mu
-\int_{SM}\divm((\mathbf X  u)\nabla^{\cdot} u)\,d\,\mu
-\int_{SM}\langle \tilde {\mathbf R}_y(\nabla^{\cdot} u),\nabla^{\cdot} u\rangle\,d\mu
\\
+\int_{SM}\Big\{\langle Y(\nabla^{\cdot}u),\nabla^{:} u\rangle
-\langle Y(y),\nabla^{\cdot}u\rangle^2
+\mathbf J((\mathbf X u)\nabla^{\cdot} u)\Big\}\,d\mu
\\
=\int_{SM}\Big\{|\nabla^{\cdot} (\mathbf X  u)|^2
-n(\mathbf X u)^2\Big\}\,d\mu.
\end{multline}

Since
$$
\divm U=U^i_{:i}=U^i_{|i}+Y^j_iU^i_{\cdot j}
=\divh U+Y^j_i U^i_{\cdot j},
$$
we have
\begin{align*}
\divm((\mathbf X  u)\nabla^{\cdot} u)
&=\divh(\mathbf X  u)\nabla^{\cdot} u)
+Y^j_i((\mathbf X u)_{\cdot j} g^{ik}u_{\cdot k}
+(\mathbf X u)g^{ik}_{\cdot j}u_{\cdot k}
+(\mathbf X u) g^{ik} u_{\cdot k\cdot j})
\\
&=\divh(\mathbf X  u)\nabla^{\cdot} u)
+\langle \nabla^{\cdot}(\mathbf Xu),Y(\nabla^{\cdot}u)\rangle,
\end{align*}
because by the symmetry argument
$$
Y^j_ig^{ik}_{\cdot j}
=-2Y^j_i g^{il}C_{lmj}g^{km}=0
$$
and
$$
Y^j_i g^{ik}u_{\cdot k\cdot j}=0.
$$

Using also (\ref{divh}), we hence have
$$
\int_{SM}\divm((\mathbf X  u)\nabla^{\cdot} u)\,d\mu
=\int_{SM}\big\{\mathbf J((\mathbf X u)\nabla^{\cdot} u)
+\langle\nabla^{\cdot } (\mathbf X u),Y(\nabla^{\cdot} u)\rangle\big\}\,d\mu.
$$

Next,
\begin{align*}
\langle \tilde {\mathbf R}_y(\nabla^{\cdot} u),\nabla^{\cdot} u\rangle
&=\big\{R^i_{kl}+(Y^i_{k|l}-Y^i_{l|k})
-(P^i_{km}Y^m_l-P^i_{lm}Y^m_k)
\\
&\quad +(Y^j_kY^i_{l\cdot j}-Y^j_lY^i_{k\cdot j})
+y_s(Y^s_l Y^i_k-Y^s_k Y^i_l)\big\}y^l u^{\cdot k}u_{\cdot i}.
\end{align*}

Now,
$$
R^i_{kl}y^l u^{\cdot k}u_{\cdot i}
=\langle \mathbf R_y(\nabla^{\cdot} u),\nabla^{\cdot} u\rangle,
$$
\begin{align*}
(Y^i_{k|l}-Y^i_{l|k})y^l u^{\cdot k}u_{\cdot i}
&=\langle (\nabla_{|y} Y)(\nabla^{\cdot} u),\nabla^{\cdot}u\rangle
-\langle\nabla_{|(\nabla^{\cdot}u)}Y(y),\nabla^{\cdot}u\rangle
\\
&=-\langle\nabla_{|(\nabla^{\cdot}u)}Y(y),\nabla^{\cdot}u\rangle
\end{align*}
by skew symmetry of $Y$ and parallelism of the fundamental tensor
with respect to~$\nabla_{|}$,
$$
(P^i_{km}Y^m_l-P^i_{lm}Y^m_k)y^l u^{\cdot k}u_{\cdot i}
=P^i_{km}Y^m_ly^lu^{\cdot k}u_{\cdot i}
=-L(Y(y),\nabla^{\cdot} u,\nabla^{\cdot}u)
$$
in view of (\ref{yp}) and (\ref{l-p}),
$$
(Y^j_kY^i_{l\cdot j}-Y^j_lY^i_{k\cdot j})y^l u^{\cdot k}u_{\cdot i}
=-2Y^j_kY^r_lg^{is}C_{srj}y^lu^{\cdot k}u_{\cdot i}
+2Y^j_lY^r_kg^{is}C_{srj}y^lu^{\cdot k}u_{\cdot i}=0
$$
by the symmetry of $C$, and
\begin{align*}
y_s(Y^s_l Y^i_k-Y^s_k Y^i_l)y^l u^{\cdot k}u_{\cdot i}
&=\langle Y(y),y\rangle\langle Y(\nabla^{\cdot} u),\nabla^{\cdot}u\rangle
-\langle Y(\nabla^{\cdot}u),y\rangle\langle Y(y),\nabla^{\cdot}u\rangle
\\
&=\langle Y(y),\nabla^{\cdot}u\rangle^2
\end{align*}
again by the skew symmetry of $Y$.

Now, (\ref{integral}) takes the form of equation (\ref{pestov-integral-final}) in the
Introduction. That is, we have proved:

\begin{Theorem}
\begin{multline}\label{pestov-integral-final-thm}
\int_{SM}\big\{|\mathbf X (\nabla^{\cdot} u)|^2
-\langle \mathbf R_y(\nabla^{\cdot} u),\nabla^{\cdot} u\rangle
-L(Y(y),\nabla^{\cdot} u,\nabla^{\cdot}u)
-\langle\nabla^{\cdot } (\mathbf X u),Y(\nabla^{\cdot} u)\rangle
\\
-2\langle Y(y),\nabla^{\cdot}u\rangle^2
+\langle \nabla^{:} u,Y(\nabla^{\cdot}u)\rangle
+\langle\nabla_{|(\nabla^{\cdot}u)}Y(y),\nabla^{\cdot}u\rangle\big\}\,d\mu
\\
=\int_{SM}\big\{|\nabla^{\cdot} (\mathbf X  u)|^2
-n(\mathbf X u)^2\big\}\,d\mu.
\end{multline}
\end{Theorem}

\begin{Remark}{\rm The identity (\ref{pestov-integral-final-thm}) is exactly identity (\ref{id})
when $n=2$. If $\phi\in C^{\infty}(TM\setminus\{0\})$ is
homogeneous of degree 1 and $n=2$, then chasing definitions we
have:
\[|\nabla^{\cdot}\phi|^2=\phi^{2}+(V\phi)^2.\]
Thus the right hand side of (\ref{pestov-integral-final-thm}) becomes
\[\int_{SM}\big\{|\nabla^{\cdot} (\mathbf X  u)|^2
-2(\mathbf X u)^2\big\}\,d\mu=
\int_{SM}\{(\G u)^2+(V\G u)^2\}\,d\mu-2\int_{SM}(\G u)^2\,d\mu\]
which is exactly the right hand side of (\ref{id}). We leave to the keen reader the task of fully verifying
that the left hand sides also coincide. When the Finsler metric is Riemannian (i.e. $I=J=0$)
and $n=2$ it is quite easy to check that (for points in $SM$):
\begin{align*}
&|\mathbf X (\nabla^{\cdot} u)|^2=(\G Vu)^2+\la^2(Vu)^2,\\
&\langle \mathbf R_y(\nabla^{\cdot} u),\nabla^{\cdot} u\rangle=(Vu)^2 K,\\
&\langle\nabla^{\cdot } (\mathbf X u),Y(\nabla^{\cdot} u)\rangle=-\la\G u\cdot Vu,\\
&\langle Y(y),\nabla^{\cdot}u\rangle=\la Vu,\\
&\langle \nabla^{:} u,Y(\nabla^{\cdot}u)\rangle=-\la \G u\cdot Vu\\
&\langle\nabla_{|(\nabla^{\cdot}u)}Y(y),\nabla^{\cdot}u\rangle=(Vu)^2 H(\la).\\
\end{align*}
Inserting these relations into the left hand side of (\ref{pestov-integral-final-thm})
we see that we get exactly the left hand side of (\ref{id}).

}
\end{Remark}

\subsection{Jacobi equation}
Let us derive a Jacobi equation. The calculations below
mimic those in the proof of \cite[Lemma 6.1.1]{She}.

Let $\phi_t:TM\setminus\{0\}\to TM\setminus\{0\}$
be the magnetic flow.
Take a curve $Z:(-\varepsilon,\varepsilon)\to TM\setminus\{0\}$
with $Z(0)=v$ and $Z'(0)=\xi$, and consider the variation
$H(s,t)=\pi(\phi_t(Z(s)))$. Set
$$
T=\de Ht,\quad U=\de Hs.
$$
Each $c_s(t)=H(s,t)$ is a magnetic geodesic; therefore,
$$
\frac{\partial^2 H^i}{\partial t^2}
+2G^i\Big(\de{H}{t}\Big)=Y^i_j\Big(\de{H}{t}\Big)\de{H^j}{t},
$$
or
\begin{equation}\label{geo}
\de{T^i}t+2G^i(T)=Y^i_j(T)T^j.
\end{equation}

Since
$$
\de{T^i}s=\de{}s\Big(\de{H^i}t\Big)
=\de{}t\Big(\de{H^i}s\Big)=\de{U^i}t,
$$
differentiating \eqref{geo} with respect to $s$
yields
\begin{align*}
\frac{\partial^2 U^i}{\partial t^2}
&=-2U^k\de{G^i}{x^k}(T)-2\de{U^l}t\de{G^i}{y^l}(T)
\\
&\quad
+\Big(U^k\de{Y^i_j}{x^k}(T)+\de{U^l}t\de{Y^i_j}{y^l}(T)\Big)T^j
+Y^i_j(T)\de{U^j}t.
\end{align*}

Note that
\begin{align*}
\de{}s\big[G^i(T)\big]
&=U^k\de{G^i}{x^k}(T)+\de{U^l}t\de{G^i}{y^l}(T),
\\
\de{}t\left[\de{G^i}{y^l}(T)\right]
&=T^k\frac{\partial^2G^i}{\partial x^k\partial y^l}
+\de{T^k}t\frac{\partial^2G^i}{\partial y^ly^k}(T)
\\
&=T^k\frac{\partial^2G^i}{\partial x^k\partial y^l}
+\big(-2G^k(T)+Y^k_m(T)T^m\big)\frac{\partial^2G^i}{\partial y^l\partial y^k}(T).
\end{align*}

Hence,
\begin{align*}
D_TD_T(U^i)&=D_T\left(\de{U^i}t+U^l\de{G^i}{y^l}(T)\right)
\\
&=\de{}t\left(\de{U^i}t+U^l\de{G^i}{y^l}(T)\right)
+\left(\de{U^k}t+U^l\de{G^k}{y^l}(T)\right)\de{G^i}{y^k}(T)
\\
&=\frac{\partial^2 U^i}{\partial t^2}
+\de{U^l}t\de{G^i}{y^l}
+U^l\de{}t\left[\de{G^i}{y^l}\right]
+\de{U^k}t\de{G^i}{y^k}
+U^l\de{G^k}{y^l}\de{G^i}{y^k}
\\
&=-2U^k\de{G^i}{x^k}-2\de{U^l}t\de{G^i}{y^l}
+\Big(U^k\de{Y^i_j}{x^k}+\de{U^l}t\de{Y^i_j}{y^l}\Big)T^j
+Y^i_j\de{U^j}t
\\
&\quad +\de{U^l}t\de{G^i}{y^l}
+U^l\left[T^k\frac{\partial^2G^i}{\partial x^k\partial y^l}
+\big(-2G^k+Y^k_mT^m\big)\frac{\partial^2G^i}{\partial y^l\partial y^k}\right]
\\
&\quad +\de{U^k}t\de{G^i}{y^k}
+U^l\de{G^k}{y^l}\de{G^i}{y^k}
\\
&=-U^k\left(2\de{G^i}{x^k}-T^j\frac{\partial^2G^i}{\partial x^j\partial y_k}
+2G^j\frac{\partial^2G^i}{\partial y^jy^k}-\de{G^i}{y^j}\de{G^j}{y^k}\right)
\\
&\quad+\Big(U^k\de{Y^i_j}{x^k}+\de{U^l}t\de{Y^i_j}{y^l}\Big)T^j
+Y^i_j\de{U^j}t
+U^lY^k_mT^m\frac{\partial^2G^i}{\partial y^l\partial y^k}.
\end{align*}

Using the identities
$$
R^i_k(T)=2\de{G^i}{x^k}-T^j\frac{\partial^2G^i}{\partial x^j\partial y_k}
+2G^j\frac{\partial^2G^i}{\partial y^jy^k}-\de{G^i}{y^j}\de{G^j}{y^k},
$$
$$
\de{U^i}t=D_T U^i-N^i_l U^l,
$$
$$
\de {Y^i_j}{x^k}=Y^i_{j|k}
+N^p_k\de{Y^i_j}{y^p}-\Gamma^i_{kp}Y^p_j+\Gamma^p_{kj}Y^i_p,
$$
$$
\frac{\partial^2G^i}{\partial y^l\partial y^k}=\Gamma^i_{jk}+L^i_{jk},
$$
$$
Y(D_T U)=Y^i_j\de{U^j}t+Y^i_j\Gamma^j_{kl}T^lU^k,
$$
we find that
$$
D_TD_T(U)=-\mathbf R_T(U)
+ Y(D_T U)+(\nabla_{|U}Y)(T)
+(\nabla_{\cdot D_T U}Y)(T)+\mathbf L(U,Y(T)),
$$
which is the Jacobi equation for the magnetic flow of a Finsler
metric. Here $\mathbf L(U,V)$ is defined by $\langle\mathbf
L(U,V),W\rangle=L(U,V,W)$.

\subsection{Index form}
Let $\gamma$ be a closed unit speed magnetic geodesic.
Let ${\mathcal A}$ and ${\mathcal C}$ be the operators on smooth vector
fields along $\gamma$ defined by
\begin{align}
{\mathcal A}(Z)&=\ddot{Z}+\mathbf R_{\dot{\gamma}}(Z)
-Y(\dot{Z})-(\nabla_{|Z}Y)(\dot{\gamma}) -(\nabla_{\cdot \dot
Z}Y)(\dot\gamma)-\mathbf L(Z,Y(\dot\gamma))   \label{A}
\\
&=\ddot Z+\mathcal C(Z) -(\nabla_{\cdot \dot
Z}Y)(\dot\gamma)-\mathbf L(Z,Y(\dot\gamma)),\nonumber
\end{align}
where
\begin{equation}
{\mathcal C}(Z):=\mathbf R_{\dot{\gamma}}(Z)
-Y(\dot{Z})-(\nabla_{|Z}Y)(\dot{\gamma}).   \label{C}
\end{equation}
If $J$ is a magnetic Jacobi field, then
\begin{equation}
{\mathcal A}(J)=0.  \label{jacobi}
\end{equation}

Let $\Lambda$ denote the $\mathbb R$-vector space of smooth vector fields
$Z:[0,T]\rightarrow TM$ along $\gamma$, such that $Z(0)=Z(T)$ and $\dot{Z}(0)=\dot{Z}(T)$.
Let $\I$ denote the quadratic form $\I:\Lambda\rightarrow \mathbb R$ defined by
\begin{equation}
\I(Z,Z)=-\int_{0}^{T}\{\left\langle {\mathcal A}(Z),Z\right\rangle+\langle Y(\dot{\gamma}),Z\rangle^{2}\}\,dt.  \label{G}
\end{equation}

Observe that
\begin{equation}\label{izz}
\I(Z,Z)=\int_0^T \big\{|\dot Z|^2-\langle \mathcal C(Z),Z\rangle
-L(Y(\dot\gamma),Z,Z)
-\langle Y(\dot\gamma),Z\rangle^2 \big\}\,dt.
\end{equation}

Indeed,
\begin{align*}
\mathbf X(\langle U,V\rangle)
&=y^k(g_{ij}U^iV^j)_{:k}
=y^k(g_{ij:k}U^iV^j+g_{ij}U^i_{:k}V^j+g_{ij}U^iV^j_{:k})
\\
&=-y^kg_{sj}Y^s_{k\cdot i}U^iV^j+\langle \mathbf X U, V\rangle
+\langle U,\mathbf X V\rangle
\\
&=-\langle (\nabla_{\cdot U}Y)(y),V\rangle+\langle \mathbf X U, V\rangle
+\langle U,\mathbf X V\rangle,
\end{align*}
where we have used the equality $g_{ij:k}=-g_{sj}Y^s_{k\cdot i}$
following from \eqref{gijk} and \eqref{yijk}.

This implies
$$
\langle \ddot Z,Z\rangle=D_{\dot\gamma}(\langle \dot Z,Z\rangle)
-|\dot Z|^2
+\langle (\nabla_{\cdot \dot Z}Y)(\dot\gamma),Z\rangle,
$$
whence \eqref{izz} is straightforward.

\begin{Lemma}[Index Lemma] Suppose the magnetic flow $\phi_{t}$ is Anosov and
let $\gamma$ be a closed magnetic geodesic with period $T$. If $Z$ is orthogonal to $\dot{\gamma}$, then
\[\I(Z,Z)\geq 0,\]
with equality if and only if $Z$ vanishes. \label{crux2}
\end{Lemma}

\begin{proof}
Let $E$ denote the weak stable or unstable subbundle of $\phi_{t}$.
It is well known (cf. \cite{PP,P3}, see \cite{CGIP} for a proof using the asymptotic Maslov
index) that the following transversality property holds:
\[E(v)\cap \mbox{\rm Ker}\,d_{v}\pi=\{0\},\]
for every $v\in SM$, where $\pi:SM \to M$ is the canonical projection.
Consider the splitting into horizontal and vertical subbundles described
in Subsection \ref{basicf}. With respect to this splitting the transversality property
can be restated as follows:
for each $v\in SM$, there exists a map $S_{v}:T_{\pi(v)}M\rightarrow T_{\pi(v)}M$ so that its
 graph is $E(v)$; moreover the correspondence $v\rightarrow S_{v}$ is continuous.

If $\xi\in E(v)$, then $J_{\xi}(t)=d\pi\circ d\phi_{t}(\xi)$ satisfies the Jacobi equation (\ref{jacobi}).
Since for all $t\in \mathbb R$,
\[\left. d\pi_{\dot{\gamma}(t)}\right|_{E(\dot{\gamma}(t))}:E(\dot{\gamma}(t))\rightarrow T_{\dot{\gamma}(t)}M\]
is an isomorphism, there exists a basis $\{\xi_{1},\dots,\xi_{n}\}$ of $E(v)$ such
that $\{J_{\xi_{1}}(t),\dots,J_{\xi_{n}}(t)\}$ is a basis of $T_{\dot{\gamma}(t)}M$ for all $t\in \mathbb R$.
Without loss of generality we may assume that $\xi_{1}=(v,S(v))$ and $J_{\xi_{1}}=\dot{\gamma}$.

Let us set for brevity $J_{i}=J_{\xi_{i}}$. Then if $Z$ is an element of $\Lambda$ we can write
\[ Z(t)=\sum_{i=1}^{n}f_{i}(t)J_{i}(t),\]
for some smooth functions $f_{1},\dots,f_{n}$ and thus,
\begin{equation}
\I(Z,Z)=-\sum_{i,j}\int_{0}^{T}\left\langle{\mathcal A}(f_{i}J_{i}),f_{j}J_{j}\right\rangle\,dt
-\int_{0}^{T}\langle Y(\dot{\gamma}),Z\rangle^{2}\,dt.      \label{bilinear}
\end{equation}
An easy computation shows that
\[{\mathcal A}(f_{i}J_{i})=\ddot{f}_{i}J_{i}
+2\dot{f}_{i}\dot{J}_{i}-\dot{f}_{i}Y(J_{i})
-\dot f_i(\nabla_{\cdot J_i}Y)(\dot\gamma)
+f_{i}{\mathcal A}(J_{i}).\]

Indeed,
$$
D_{\dot\gamma}D_{\dot\gamma}(f_iJ_i)
=\ddot f_iJ_i+2\dot f_i\dot J_i+f_i\ddot J_i,
$$
$$
\mathbf R_{\dot\gamma}(f_iJ_i)=f_i\mathbf R_{\dot\gamma}(J_i),
$$
$$
Y(D_{\dot\gamma}(f_iJ_i))
=\dot f_i Y(J_i)+f_i Y(\dot J_i),
$$
$$
(\nabla_{|f_iJ_i}Y)(\dot\gamma)=f_i(\nabla_{|J_i}Y)(\dot\gamma),
$$
$$
(\nabla_{\cdot D_{\dot\gamma}(f_iJ_i)}Y)(\dot\gamma)
=\dot f_i(\nabla_{\cdot J_i}Y)(\dot\gamma)
+f_i(\nabla_{\cdot \dot J_i}Y)(\dot\gamma),
$$
$$
\mathbf L(f_iJ_i,Y(\dot\gamma))=f_i\mathbf L(J_i,Y(\dot\gamma)).
$$

Since $J_{i}$ satisfies equation (\ref{jacobi}),
it follows that ${\mathcal A}(J_{i})=0$ and hence
\[\langle{\mathcal A}(f_{i}J_{i}),J_{j}\rangle
=\ddot{f}_{i}\langle J_{i},J_{j}\rangle
+2\dot{f}_{i}\langle\dot{J}_{i},J_{j}\rangle
-\dot{f}_{i}\langle Y(J_{i}),J_{j}\rangle
-\dot f_i\langle(\nabla_{\cdot J_i}Y)(\dot\gamma),J_j\rangle.\]

Observe that since $E$ is a Lagrangian subspace,
\[\langle J_{i},\dot{J}_{j}\rangle
-\langle\dot{J}_{i},J_{j}\rangle
+\langle Y(J_{i}),J_{j}\rangle=0,\]
and then
\[\langle{\mathcal A}(f_{i}J_{i}),J_{j}\rangle
=\frac{d}{dt}(\dot{f}_{i}\langle J_{i},J_{j}\rangle).\]
Now we can write
\[\int_{0}^{T}\langle{\mathcal A}(f_{i}J_{i}),f_{j}J_{j}\rangle\,dt
=\left.\langle\dot{f}_{i}J_{i},f_{j}J_{j}\rangle\right|_{0}^{T}
-\int_{0}^{T}\langle\dot{f}_{i}J_{i},\dot{f}_{j}J_{j}\rangle\,dt.\]
Combining the last equality with (\ref{bilinear}) we obtain
\[\I(Z,Z)=\int_{0}^{T}\bigg|\sum_{i=1}^{n}\dot{f}_{i}J_{i}\bigg|^{2}\,dt
-\bigg\langle\sum_{i=1}^{n}\dot{f}_{i}J_{i},Z\bigg\rangle\bigg|_{0}^{T}
-\int_{0}^{T}\langle Y(\dot{\gamma}),Z\rangle^{2}\,dt.\]
But $\dot{Z}(0)=\dot{Z}(T)$ and $\dot{Z}=\sum_{i=1}^{n}\dot{f}_{i}J_{i}+\sum_{i=1}^{n}f_{i}\dot{J}_{i}$, therefore
\[\bigg\langle\sum_{i=1}^{n}\dot{f}_{i}J_{i},Z\bigg\rangle\bigg|_{0}^{T}
=-\bigg\langle\sum_{i=1}^{n}f_{i}\dot{J}_{i},Z\bigg\rangle\bigg|_{0}^{T}.\]
Note that $\dot{J}_{i}(t)=S_{\dot{\gamma}(t)}J_{i}(t)$, hence
\[\sum_{i=1}^{n}f_{i}\dot{J}_{i}=S\left(\sum_{i=1}^{n}f_{i}J_{i}\right)=S(Z),\]
which implies
\[\bigg\langle\sum_{i=1}^{n}f_{i}\dot{J}_{i},Z\bigg\rangle\bigg|_{0}^{T}
=\bigg\langle S(Z),Z\bigg\rangle \bigg|_{0}^{T}=0.\]
Then
\begin{equation}
\I(V,V)=\int_{0}^{T}\bigg| \sum_{i=1}^{n}\dot{f}_{i}J_{i}\bigg| ^{2}\,dt
-\int_{0}^{T}\langle Y(\dot{\gamma}),Z\rangle^{2}\,dt.
\label{finvi}
\end{equation}
Now let
\[W:=\sum_{i=2}^{n}\dot{f}_{i}J_{i}.\]
Since $J_1=\dot{\gamma}$ we have:
\[\bigg\langle
\sum_{i=1}^{n}\dot{f}_{i}J_{i},\sum_{i=1}^{n}\dot{f}_{i}J_{i}\bigg\rangle
=\langle \dot{f}_{1}\dot{\gamma}+W, \dot{f}_{1}\dot{\gamma}+W\rangle
=\dot{f}_{1}^{2}+2\dot{f}_{1}\langle \dot{\gamma},W\rangle+\langle W,W\rangle.\]
Differentiating $\langle Z,\dot{\gamma}\rangle=0$ we get
\[\langle \dot{Z},\dot{\gamma}\rangle+\langle Z,Y(\dot{\gamma})\rangle=0.\]
But
\[\langle \dot{Z},\dot{\gamma}\rangle=\left\langle \sum_{i=1}^{n}\dot{f}_{i}J_{i},\dot{\gamma}\right\rangle
=\dot{f}_{1}+\langle W,\dot{\gamma}\rangle\]
since $\langle \dot{J}_{i},\dot{\gamma}\rangle=0$ for all $i$. Therefore
\[\langle Y(\dot{\gamma}),Z\rangle^{2}=\dot{f}_{1}^{2}+2\dot{f}_{1}\langle W,\dot{\gamma}\rangle+
\langle W,\dot{\gamma}\rangle^{2}.\]
Thus
\[\left\langle \sum_{i=1}^{n}\dot{f}_{i}J_{i},\sum_{i=1}^{n}\dot{f}_{i}J_{i}\right\rangle
-\langle Y(\dot{\gamma}),Z\rangle^{2}=\langle W,W\rangle-\langle W,\dot{\gamma}\rangle^{2}.\]
If we let $W^{\perp}$ be the orthogonal projection of $W$ to $\dot{\gamma}^{\perp}$, the last equation
and (\ref{finvi}) give:
\[\I(Z,Z)=\int_{0}^{T}\| W^{\perp}\|^{2}\,dt\geq 0\]
with equality if and only if $W^{\perp}$ vanishes identically. But if $W^{\perp}$ vanishes, then
\[-\langle W,\dot{\gamma}\rangle\dot{\gamma}+\sum_{i=2}^{n}\dot{f}_{i}J_{i}=0\]
which implies that the functions $f_{i}$ are constant for $i\geq 2$. Thus
$Z$ is of the form $f_{1}\dot{\gamma}+ J$ where $J$ is a magnetic Jacobi field.
If we let $J^{\perp}$ be the orthogonal projection of $J$ to $\dot{\gamma}^{\perp}$, then
$Z=J^{\perp}$.
Now write
\[J=x\dot{\gamma}+J^{\perp}\]
A simple calculation shows that ${\mathcal
A}(x\dot{\gamma})=D_{\dot\gamma}(\dot{x}\dot{\gamma})$ with
$\dot{x}=\langle J,Y(\dot{\gamma})\rangle=\langle
J^{\perp},Y(\dot{\gamma})\rangle$. Hence
\[0={\mathcal A}(J)={\mathcal A}(J^{\perp})+D_{\dot\gamma}\left(\langle J^{\perp},Y(\dot{\gamma})\rangle\dot{\gamma}\right).\]
The fact that $J^{\perp}$ satisfies this second order differential equation together with
$J^{\perp}(0)=J^{\perp}(0)$ and $\dot{J^{\perp}}(T)=\dot{J^{\perp}}(T)$ implies that
$J^{\perp}$ is periodic with period $T$. Hence
$\dot{x}$ is also a periodic function of period $T$ which implies that $\| J\| $ grows at most
linearly with $t$. However, since the closed orbits of $\phi_{t}$ are hyperbolic
the only Jacobi fields with that type of growth are those given by constant multiples of $\dot{\gamma}$.
Since $Z$ is orthogonal to $\dot{\gamma}$, $Z$ must vanish.
\end{proof}

\subsection{End of the proof of Theorem B}
Define
$$
\tilde{\mathcal C}(V)=\mathbf R_y(V)-Y(\mathbf X V)-(\nabla_{|V} Y)(y).
$$

Then the following holds:
\begin{align*}
\langle \tilde{\mathcal C}(\nabla^{\cdot}u),\nabla^{\cdot}u\rangle
&=\langle \mathbf R_y(\nabla^{\cdot} u),\nabla^{\cdot} u\rangle
+\langle \mathbf X(\nabla^{\cdot}u),Y(\nabla^{\cdot} u)\rangle
-\langle(\nabla_{|(\nabla^{\cdot} u)}Y)(y),\nabla^{\cdot}u\rangle
\\
&=\langle \mathbf R_y(\nabla^{\cdot}u),\nabla^{\cdot} u\rangle
+\langle \nabla^{\cdot}(\mathbf X  u)-\nabla^{:}u
-\langle Y(y),\nabla^{\cdot}u\rangle y,
Y(\nabla^{\cdot}u)\rangle
\\
&\quad-\langle(\nabla_{(\nabla^{\cdot} u)}Y)(y),\nabla^{\cdot}u\rangle
\\
&=\langle \mathbf R_y(\nabla^{\cdot}u),\nabla^{\cdot} u\rangle
+\langle \nabla^{\cdot}(\mathbf X  u),
Y(\nabla^{\cdot}u)\rangle
-\langle\nabla^{:}u,Y(\nabla^{\cdot}u)\rangle
\\
&\quad+\langle Y(y),\nabla^{\cdot}u\rangle^2
-\langle(\nabla_{|(\nabla^{\cdot} u)}Y)(y),\nabla^{\cdot}u\rangle.
\end{align*}

Suppose $\G u=h\circ\pi+\theta$. From
(\ref{pestov-integral-final-thm}) and Lemma \ref{int-nabla}
we infer that
\begin{equation}
\int_{SM}\big\{|\mathbf X \nabla^{\cdot} u|^2
-\langle \tilde{\mathcal C}( \nabla^{\cdot} u), \nabla^{\cdot} u\rangle
-L(Y(y), \nabla^{\cdot} u, \nabla^{\cdot} u)
-\langle Y(y), \nabla^{\cdot} u\rangle^2\big\}\,d\mu\leq 0.
\label{menqc}
\end{equation}

Given a closed unit-speed magnetic geodesic $\gamma:[0,T]\to M$
consider the smooth vector field $Z:[0,T]\to TM$ along $\gamma$ given by
$Z:= \nabla^{\cdot} u(\ga,\dot{\ga})$. Note that $Z$ is orthogonal to
$\dot{\ga}$ because $u$ is homogeneous of degree zero.

The Index Lemma \ref{crux2} tells us that
\begin{equation}\label{in_gN}
\int_0^T \big\{|\dot Z|^2-\langle \mathcal C(Z),Z\rangle
-L(Y(\dot\gamma),Z,Z)
-\langle Y(\dot\gamma),Z\rangle^2\big\} \,dt\geq 0
\end{equation}
for every closed magnetic geodesic $\ga$.

Since the flow is Anosov, the invariant measures supported on
closed orbits are dense in the space of all invariant measures
on $SM$.
Therefore, the above yields
$$
\int_{SM}\big\{|\mathbf X \nabla^{\cdot} u|^2
-\langle \tilde{\mathcal C}( \nabla^{\cdot} u), \nabla^{\cdot} u\rangle
-L(Y(y), \nabla^{\cdot} u, \nabla^{\cdot} u)
-\langle Y(y), \nabla^{\cdot} u\rangle^2\big\}\,d\mu\geq 0.
$$
Combining this (\ref{menqc}), we find that
\begin{equation}
\int_{SM}\big\{|\mathbf X \nabla^{\cdot} u|^2
-\langle \tilde{\mathcal C}( \nabla^{\cdot} u), \nabla^{\cdot} u\rangle
-L(Y(y), \nabla^{\cdot} u, \nabla^{\cdot} u)
-\langle Y(y), \nabla^{\cdot} u\rangle^2\big\}\,d\mu= 0.
\label{in_smN}
\end{equation}

By the non-negative version of the Liv\v sic theorem, proved independently by M. Pollicott and R. Sharp and by
A. Lopes and P. Thieullen (see \cite{LT, PS}),
we conclude from (\ref{in_gN}) and (\ref{in_smN}) that

\[\int_0^T \big\{|\dot Z|^2-\langle \mathcal C(Z),Z\rangle
-L(Y(\dot\gamma),Z,Z)
-\langle Y(\dot\gamma),Z\rangle^2\big\} \,dt= 0\]
for every closed magnetic geodesic $\gamma$.
Applying again the Index Lemma \ref{crux2}, we see that
$\nabla^{\cdot} u$ vanishes on all closed magnetic geodesics.
Since the latter are dense in $SM$, the function $\nabla^{\cdot} u$
vanishes on all of $SM$.
This means that $u=f\circ\pi$ where $f$ is a smooth function on $M$. But in this case, since
$d\pi_{(x,v)}(\G)=v$ we have
$\G(u)=df_{x}(v)$ and Theorem B follows.

\section{Proof of Theorem C}

Suppose the magnetic flow $\phi$ of the pair $(F,\Omega)$ has an Anosov splitting
$$E^s\oplus E^u\oplus \re\G$$
 of class $C^1$ and suppose also that $\Omega$ is exact.
Let $\tau$ denote the one-form that vanishes
 on $E^{s}\oplus E^{u}$ and takes the value one on the vector
 field $\G$. If the splitting is of class $C^{1}$ then
 $\tau$ is also of class $C^{1}$ and $d\tau$ is a continuous
 2-form invariant under the magnetic flow.
 U. Hamenst\"adt showed in \cite{Ha}, for the geodesic flow case, that any continuous
 invariant exact 2-form must be a constant multiple of the symplectic form provided
 that the splitting is of class $C^{1}$. Hamenst\"adt's proof carries over to the case of
magnetic flows {\it without} major changes, provided that $\Omega$ is an exact
form $d\theta$ (see the appendix of \cite{Pa}). Recall from the introduction that
the symplectic form on $TM\setminus \{0\}$
is given by $\omega_{0}+\pi^{*}\Omega$, where
$\omega_{0}=\ell_{F}^{*}(-d\lambda)$ ($\ell_{F}$ is the Legendre transform of $F^2/2$ and $\la$ is
the Liouville 1-form of $T^*M$). It follows that there exists
 a constant $c$ such that:
 \[d\tau=c(\omega_{0}+\pi^{*}\Omega),\]
 and thus
 \[d(\tau+c\ell_{F}^{*}\lambda-c\pi^{*}\theta)=0.\]
Let us write
\[\varphi:=\tau+c\ell_{F}^{*}\lambda-c\pi^{*}\theta.\]
Then $\varphi$ is a smooth {\it closed} 1-form. Since on $SM$
$\ell_{F}^{*}\lambda(\G)=1$\footnote{Using the expressions in
Subsection \ref{idloc} we see that
$\ell_{F}^{*}\lambda(\G)=g_{ij}y^iy^j$.} we obtain
\begin{equation}
\varphi(\G)(x,v)=1+c-c\theta_{x}(v). \label{p}
\end{equation}
It is well known that the map $\pi^{*}:H^{1}(M,\re)\to H^{1}(SM,\re)$
is an isomorphism (provided that $M$ is not
diffeomorphic to a 2-torus).
 Therefore there exist a {\it closed} smooth 1-form $\delta$
 in $M$
and a smooth function $u:SM\to \re$ such that
\[\varphi=\pi^{*}\delta+du.\]
Hence equation (\ref{p}) gives:
\begin{equation}\label{p1}
\G(u)+\delta_{x}(v)=1+c-c\theta_{x}(v).
\end{equation}
Integrating the last equality with respect to the (normalized) Liouville measure $\mu$ and using that
the magnetic flow leaves $\mu$ invariant we have
\[0=1+c-c\int_{SM}\theta\,d\mu-\int_{SM}\delta\,d\mu.\]
By Lemma \ref{int-nabla}
\[\int_{SM}\theta\,d\mu=\int_{SM}\delta\,d\mu=0\]
and thus $c=-1$. Replacing in (\ref{p1}) we finally obtain
\begin{equation}
\delta_{x}(v)+\G(u)(x,v)=\theta_{x}(v).  \label{ch}
\end{equation}
We can now apply Theorem B to conclude that $\theta$ is a closed
form, i.e., $\Omega$ vanishes identically.\footnote{Alternatively,
we could have applied Theorem B directly to equation (\ref{p1}) to
conclude that $c=-1$ and $\theta$ is exact.}

\end{document}